\documentclass{amsart}
\usepackage{etex}
\usepackage{graphicx}
\usepackage{amssymb}
\usepackage{epstopdf}
\usepackage{nicefrac}
\usepackage{enumerate}
\usepackage{tabularx}
\usepackage{hyperref}
\usepackage{pst-all}
\usepackage{tikz}
\usepackage{xcolor}

\usetikzlibrary{trees}

\newtheorem{thm}{THEOREM}[section]

\newtheorem{cor}[thm]{COROLLARY}
\newtheorem{defn}[thm]{DEFINITION}

\newtheorem{lemma}[thm]{LEMMA}

\newtheorem{prop}[thm]{PROPOSITION}

\newtheorem{remark}[thm]{REMARK}

\newcommand{\cA}{{\mathcal A}}
\newcommand{\cB}{{\mathcal B}}

\newcommand{\cF}{{\mathcal F}}
\newcommand{\cG}{{\mathcal G}}
\newcommand{\cI}{{\mathcal I}}
\newcommand{\cJ}{{\mathcal J}}
\newcommand{\cM}{{\mathcal M}}

\newcommand{\cQ}{{\mathcal Q}}
\newcommand{\cS}{{\mathcal S}}
\newcommand{\cT}{{\mathcal T}}
\newcommand{\cU}{{\mathcal U}}
\newcommand{\cV}{{\mathcal V}}
\newcommand{\cX}{{\mathcal X}}
\newcommand{\diam}{{\rm{diam}}}
\newcommand{\ds}{\displaystyle}
\newcommand{\e}{{\epsilon}}
\newcommand{\fM}{{\mathfrak{M}}}
\newcommand{\fP}{{\mathfrak{P}}}
\newcommand{\fX}{{\mathfrak{X}}}
\newcommand{\G}{\Gamma}
\newcommand{\mN}{{\mathbb N}}
\newcommand{\mR}{{\mathbb R}}
\newcommand{\mZ}{{\mathbb Z}}

\begin{document}

\title{Hausdorff dimension in graph matchbox manifolds}

\begin{abstract}{We study the Hausdorff and the box dimensions of closed invariant subsets of the space of pointed trees, equipped with a pseudogroup action.  This pseudogroup dynamical system can be regarded as a generalization of a shift space. We show that the Hausdorff dimension of this space is infinite, and the union of closed invariant subsets with dense orbit and non-equal Hausdorff and box dimensions is dense in this space.

We apply our results to the problem of embedding laminations into differentiable foliations of smooth manifolds. One of necessary conditions for the existence of such an embedding is that the lamination must admit a bi-Lipschitz embedding into a manifold. A suspension of the pseudogroup action on the space of pointed graphs gives an example where this condition is not satisfied, with Hausdorff dimension of the space of pointed trees being the obstruction to the existence of such a bi-Lipschitz embedding.

}
\end{abstract}


\author{Olga Lukina}
\email{olga.lukina@univie.ac.at}
\address{Institute for Mathematics, University of Vienna, Oskar-Morgenstern-Platz 1, 1090 Vienna, Austria}

\thanks{Version date: September 18, 2016. Revisions: November 26, 2021 and December 21, 2021.}

\date{}


\maketitle

\section{Introduction} \label{sec-intro}
 
Let $G$ be a finitely generated group, and $\Sigma:G \to \cA$ be the space of all maps from $G$ to a finite set $\cA$. There is a natural action of $G$ on $\Sigma$, and $\Sigma$ with this action is called the \emph{Bernouilli shift} \cite{CoorPapad}. Bernouilli shifts are well-known and well-studied, especially in the case $G = \mN$ or $G = \mZ$, and they find their applications in many branches of dynamical systems theory and other areas of mathematics. They have interesting properties, which make them useful in applications, such as positive entropy, positive Hausdorff dimension, and others.

Ghys \cite{Ghys1999} considered a dynamical system, where the phase space is the space of non-compact trees of bounded valence of vertices and with one distinguished vertex. Topologically the space of such pointed trees is a Cantor set. There is a partial action of a free group $F_n$, $n \geq 2$, on this space, which can be thought of as moving the distinguished vertex along the edges of trees. Because of the action, one intuitively thinks of the space of pointed trees with a partial action of a free group as a generalized Bernouilli shift. 

Given a pseudogroup dynamical system on a topological space $X$ it is natural to ask, what are the properties of dynamically defined subsets of this space? Here by `dynamically defined' we mean the subsets which are closed, invariant under the action of the pseudogroup, and which contain a dense orbit of the pseudogroup action. These subsets need not be minimal with respect to the properties of being closed and invariant, and so they form a hierarchy ordered by inclusion. This hierarchy in codimension one foliations, known as `theory of levels', was studied by Cantwell and Conlon \cite{CC1}, Hector \cite{Hector1983}, Nishimori \cite{Nishimori1977} and others, see \cite{Lukina2012} for references. 

The hierarchy of closed invariant subsets with a dense orbit in the space of pointed trees, given by inclusions, was studied by the author in \cite{Lukina2012}. Specific examples of such subsets were considered in ~Ghys \cite{Ghys1999}, Blanc \cite{Blanc2001}, ~Alcalde~Cuesta, ~Lozano~Rojo and ~Macho~Stadler \cite{ALM2009a}. Lozano Rojo \cite{Lozano2009} showed that the space of pointed trees contains closed invariant subsets with a dense orbit which are not uniquely ergodic. In Lukina \cite{Lukina2012}, it was shown that the structure of orbits of the action on the space of pointed trees is remarkably similar to the structure of the orbits in the full shift space. In Lozano Rojo and Lukina \cite{LozanoLukina2012},  it was shown that a full shift on any finitely generated group and any finite alphabet can be equivariantly embedded into the space of pointed trees. 

The next natural question to ask is, what are the dimension properties of the closed invariant subsets with a dense orbit in the space of pointed trees? For example, it is well known that such dynamically defined subsets of the shift space have equal Hausdorff and box dimensions. Since the shift space equivariantly embeds into the space of pointed trees \cite{LozanoLukina2012}, one can expect that the space of pointed trees contains subsets with equal box and Hausdorff dimensions, but these subsets are very special. In this paper we show that the space of pointed trees exhibits much more rich and diverse properties than that. In particular, we show that it has infinite Hausdorff dimension, and that any open set in this space intersects a) an orbit whose closure has zero Hausdorff dimension and positive box dimension, b) an orbit whose closure has positive Hausdorff dimension.

Our computations have the following application. One can suspend the action on the space of pointed trees, obtaining a lamination which we call the \emph{space of graph matchbox manifolds}. This lamination is an example of a lamination with bi-Lipschitz holonomy pseudogroup, which does not admit a bi-Lipschitz embedding as a foliated subset of a smooth manifold with a $C^{\infty,r}$-foliation, for $r \geq 1$, see the discussion below for details.

\bigskip
We now give precise statements of our results.

Let $F_n$ be a free group on $n$ generators, $n \geq 2$. We assume that the generating set $F_n^{1+}$ is not symmetric, that is, if $a \in F_n^{1+}$ then $a^{-1} \notin F_n^{1+}$. We denote by $\cG_n$ the Cayley graph of $F_n$, determined by the choice of generators $F_n^{1+}$. Let $X_n$ be the set of all connected, simply-connected non-compact subgraphs of $\cG_n$ containing the identity vertex. Every such graph $T \in X_n$ is a pointed metric space, with complete length metric associated to the standard length structure on $\cG_n$.

The set $X_n$ of pointed trees is given a \emph{ball metric} as follows. Given any two trees $T$ and $T'$ in $X_n$, the distance $d_n(T,T') = e^{-m}$, where $m$ is the maximal radius such that the closed balls of this radius in $T$ and $T'$ correspond to the same labelled pattern of edges.

An element $g \in F_n$ acts on a tree $T$ if $g \in V(T)$, where $V(T) \subset F_n$ denotes the set of vertices of $T$. The action takes $T$ to the tree $T'$ with exactly the same geometry, but where the distinguished vertex is positioned differently (see Section \ref{gmm-section} for details). Each element $g \in F_n$ gives rise to a homeomorphism $\gamma_g$, defined on a clopen subset of $X_n$. So one obtains a pseudogroup dynamical system $(X_n,d_n,\G_n)$, where $\G_n = \{\gamma_g\}_{g \in F_n}$. 

Let $T \in X_n$, and let $O(T)$ be the orbit of $T$ under the action of the pseudogroup $\G_n$. Then the closure $\cX = \overline{O(T)}$ is a closed subset of $X_n$. This subset $\cX$ is transitive, that is, it contains a dense orbit $O(T)$, and it is invariant under the action of $\G_n$. We are interested in the box and the Hausdorff dimensions of such subsets, in particular, if they coincide, and if they are finite or infinite.

The \emph{lower box dimension} $\underline{\dim}_B(Z)$ and the \emph{upper box dimension} $\overline{\dim}_B(Z)$ are often used to estimate the Hausdorff dimension $\dim_H(Z)$. If $\underline{\dim}_B(Z) = \overline{\dim}_B(Z)$, then we say that the box dimension $\dim_B(Z)$ is defined for $Z$. The box dimension is also called the Minkowski dimension, and, in general,
  \begin{align}\label{eq-3d}\dim_H(Z) \leq \underline{\dim}_B(Z) \leq \overline{\dim}_B(Z).\end{align}
 For the set of pointed trees, we observe the following. 

\begin{prop}\label{prop-subsets}
The space $(X_n,d_n)$ contains closed subsets, invariant under the action of the pseudogroup $\G_n$, with the following properties.
\begin{enumerate}
\item The Hausdorff dimension of $(X_n,d_n)$ is infinite,
   \begin{align*} d_H(X_n) = \infty. \end{align*} 
\item  Let $\cI $ be the set of closed transitive invariant subsets of $X_n$ with zero Hausdorff dimension and positive box dimension. Then the union $\bigcup_{\cM \in \cI} \cM$ is a dense subset of $X_n$.
\item Let $\cJ $ be the set of closed transitive invariant subsets of $X_n$ with non-zero Hausdorff dimension. Then the union $\bigcup_{\cM \in \cJ} \cM$ is a dense subset of $X_n$.
\end{enumerate}
\end{prop}

Statement $(1)$ in Proposition \ref{prop-subsets} is a consequence of the space $(X_n,d_n)$ containing a subset with a lot of regularity, which allows to compute its Hausdorff dimension using only a basic counting technique. The proofs of items $2)$ and $3)$ in Proposition \ref{prop-subsets} are constructive. We give a brief outline of the ideas involved in these proofs. The details are spelled out in Section \ref{sec-subsets}.

Given a clopen set $U$ in $X_n$, there is a finite graph $G_{U}$, such that for every $T \in U$ the pattern around the distinguished vertex is that of $G_{U}$. We build a graph $T_{U}$ by attaching to $G_U$ a specific subgraph of $F_n$ at the boundary of $G_U$. This subgraph is chosen in such a way that the closed subset $\cX_U = \overline{O(T_U)}$ has positive box dimension. By the properties of the construction, $\cX_U$ consists of only $2$ orbits. Therefore, $\cX_U$ is countable and so has zero Hausdorff dimension. Since the pattern around the distinguished vertex in $T_U$ is that of $G_U$, we have that $T_U \in U$. Since $U$ is arbitrary, it follows that the union of all closed transitive invariant subsets with zero Hausdorff and positive box dimension is dense in $X_n$.

The proof of 3) uses embeddings of shift spaces on $\mZ$ into the space of pointed trees $X_n$ developed in \cite{LozanoLukina2012}, and the fusion construction from \cite{Lukina2012}. It is well-known that the upper and the lower box dimensions and the Hausdorff dimension of a shift space are equal and positive. It is also well-known that a shift space has elements with transitive orbits. We use the construction from \cite{LozanoLukina2012} to build an embedding of a shift space into $X_n$ by associating to each element in the shift space a graph in $X_n$.  This embedding is not unique and depends on choices. The following lemma shows that such embeddings preserve the positivity of the Hausdorff dimension, and the equality of the Hausdorff and the box dimensions.

\begin{lemma}\label{lemma-holderembedding}
Let  $\cA$ be a finite alphabet, and $\Sigma = \{\sigma: \mZ \to \cA\}$ be a shift space. Let $\Phi: \Sigma \to X_n$ be an equivariant embedding as in \cite{LozanoLukina2012}. Then $\Phi$ is H\"older continuous, and 
  $$0<\dim_H \, \Phi(\Sigma) = \underline{\dim}_B \, \Phi(\Sigma) =  \overline{\dim}_B \, \Phi(\Sigma) < \infty.$$
\end{lemma}

To complete the proof of 3) we use the fusion construction from \cite{Lukina2012}. Let $T_s$ be a graph which is an image of an element of the shift space with transitive orbit under the embedding from Lemma \ref{lemma-holderembedding}. For a clopen set $U \subset X_n$, choose $T_U \in U$. The fusion is a procedure to construct a graph $T$ in $X_n$ such that $\overline{O(T_U)}, \overline{O(T_s)} \subset \overline{O(T)}$. It follows that $\overline{O(T)}$ intersects $U$ and has positive Hausdorff dimension. We note that the box dimension of $\overline{O(T)}$ need not be equal to its Hausdorff dimension. 

The image of an embedding in Lemma \ref{lemma-holderembedding} is a subset for which the Hausdorff and the box dimension are positive and equal, although we do not know if the union of subsets with this property is dense in $X_n$.

\medskip
One  application of Proposition \ref{prop-subsets} is to the problem of embedding laminations into $C^{\infty,r}$ foliations of smooth manifolds, for $r \geq 1$.

The pseudogroup dynamical system $(X_n,d_n,\G_n)$, associated to a free group $F_n$, $n \geq 2$, can be suspended to produce a lamination $\fM_n$ with leaves of dimension $2$ or higher \cite{Ghys1999,ALM2009a,Blanc2001,LozanoLukina2012,Lukina2012}.  
A \emph{graph matchbox manifold} is a compact connected subspace $\cM$ of $\fM_n$ which is a suspension of the closure $\overline{O(T)} \subset X_n$ of the orbit of a tree $T$ under the action of the pseudogroup $\G_n$. Since $X_n$ contains dense orbits of the $\G_n$-action \cite{Lukina2012}, then $\fM_n$ is a graph matchbox manifold in its own right.

\medskip
Laminations provide models for the study of attractors in dynamical systems and, more generally, minimal sets of foliations of smooth manifolds. For example, Williams \cite{Williams1967,Williams1974}
showed that the dynamics on an expanding hyperbolic attractor is conjugate to the shift map on the inverse limit of a map of some branched manifold onto itself. This inverse limit is called a \emph{generalized solenoid}, and it is an example of a lamination. Brown \cite{Brown2010} showed that the dynamics on a non-expanding topologically mixing codimension $1$ hyperbolic attractor is conjugate to an automorphism of an inverse limit of tori, which is also an example of a lamination. 

Therefore, it is natural to ask if a given lamination can be realized as a foliated subset of a $C^{\infty,r}$-foliation of a $C^{\infty}$ manifold, for $r \geq 1$.

Current literature contains just a few examples of embeddings of laminations into foliations, especially into $C^{\infty,r}$-foliations, for $r \geq 1$. Here $C^{\infty,r}$ means that the leaves of a foliation are smooth manifolds, but the manifold $M$ admits a \emph{foliated} atlas  of class   $C^r$, $r \geq 1$. We refer the reader to Section \ref{sec-application} for an overview of the literature. For graph matchbox manifolds, Lozano \cite{Lozano2009} and Hector \cite{Hector2014} considered their embeddings into $C^{\infty,0}$-foliations, i.e. the case when a foliation of a manifold is just continuous, and not differentiable. 

The requirement that the foliation is $C^{\infty,r}$, for $r \geq 1$, imposes additional conditions on a lamination. For example, the classical Denjoy foliation of a $2$-torus, see Denjoy \cite{Denjoy}, can only appear in $C^0$ or $C^1$ foliations. Further,
 McDuff \cite{McDuff1981} showed that Cantor sets with certain metrics cannot occur as minimal sets of $C^1$-homeomorphisms of a circle, and so their suspensions cannot occur as exceptional minimal sets of $C^1$-foliations of a $2$-torus. 

\medskip
The conditions imposed on a lamination and on the embedding map by the requirement that the foliation of a manifold is $C^{\infty,r}$, for $r \geq 1$, in any codimension, were studied by Hurder \cite{Hurder2013}. 

Recall, that the map $f: Y \to Z$ of metric spaces $(Y,d_Y)$ and $(Z,d_Z)$ is \emph{bi-Lipschitz}, if there exists a constant $K>0$ such that for any $y,y' \in Y$ we have
\begin{align*} \frac{1}{K} d_{Y}(y,y') \leq d_Z(f(y),f(y')) \leq K d_Y(y,y') . \end{align*}
If a lamination $\fM$ can be realized as a subset of a $C^{\infty,r}$ foliation of a manifold $M$, for $r \geq 1$, then it must satisfy the following two conditions.

\begin{enumerate}
\item \label{cons-1} $\fM$ must admit a foliated atlas $\cQ$ and a metric $d_{\fM}$, such that the generators of the holonomy pseudogroup associated to $\cQ$ are bi-Lipschitz maps with respect to the metric induced on the transversal by $d_{\fM}$.
\item  \label{cons-2}There must exist a homeomorphism onto its image $\phi: \fM \to M$ which is a bi-Lipschitz map.
\end{enumerate}

We show in Section \ref{gmm-section}, that the pseudogroup $\G_n$ acting on the space of pointed trees $(X_n,d_n)$ is generated by bi-Lipschitz maps. Then we have the following consequence of Proposition \ref{prop-subsets}.

\begin{prop}\label{prop-main3}
Let $F_n$, $n \geq 2$, be a free group on $n$ generators, and $(X_n,d_n)$ be the associated space of pointed trees with ball metric. Let $\fM_n$ be the suspension of the natural pseudogroup action on $(X_n,d_n)$ as in \cite{Ghys1999}. Then for $r\geq 1$, the lamination $\fM_n$ cannot be embedded as a foliated subset of a $C^\infty$-manifold with a $C^{\infty,r}$  foliation of any codimension by a bi-Lipschitz homeomorphism onto its image.
\end{prop}

Proposition \ref{prop-main3} is a consequence of the facts, that bi-Lipschitz embeddings preserve the Hausdorff dimension of sets \cite{Pesin1997}, and that the Hausdorff dimension of a subset of a manifold is equal to at most the dimension of the ambient manifold. Proposition \ref{prop-main3} gives an example of a lamination $\fM_n$ with metric $d_{\fM}$ which satisfies Condition \eqref{cons-1}, but does not satisfy Condition \eqref{cons-2}.

Although one can construct examples of metrics on Cantor sets, so that their Hausdorff dimension is infinite, in all examples known to the author, it seems problematic to equip these Cantor sets with a pseudogroup dynamical system so that the system is compactly generated, and the generators are bi-Lipschitz maps relative to the metric. Thus the space of graph matchbox manifolds in Proposition \ref{prop-main3}  seems to be the first known example of its kind. 

Hurder \cite{Hurder2013} uses the result of Proposition \ref{prop-main3}  as an example of a lamination which does not admit a bi-Lipschitz embedding into a $C^{\infty,r}$-foliation, where $r \geq 1$, not only for the given ball metric $d_n$, but for \emph{any} metric on $\fM_n$, such that  the holonomy pseudogroup of $\fM_n$ is generated by bi-Lipschitz maps relative to this metric. His argument is based on the fact, that if $\fM_n$ has infinite Hausdorff dimension, then it has infinite pseudogroup entropy. The fact that entropy is infinite depends neither on the metric on $\fM_n$, nor on the generating set of $\G_n$ (see Ghys, Langevin and Walczak \cite{GLW1988}). It follows that $\fM_n$ cannot be embedded into a $C^{\infty,r}$-foliation by a bi-Lipschitz homeomorphism onto its image, for $r \geq 1$, since the entropy of a differentiable foliation of a smooth finite-dimensional manifold is always finite (see Ghys, Langevin and Walczak \cite{GLW1988}). The result of Proposition \ref{prop-main3} in this paper is vital for this conclusion.

There are many open problems which concern the Hausdorff dimension of graph matchbox manifolds. Here are just a few of them:

\begin{enumerate}
\item The example in Proposition \ref{prop-main3}  of a lamination with compactly generated bi-Lipschitz holonomy pseudogroup, which does not admit a bi-Lipschitz embedding into a $C^{\infty,r}$-foliated manifold for $r \geq 1$, is not minimal. It would be interesting to find a minimal example with similar properties.
\item Find a graph $T \subset X_n$ such that the lower box dimension $\underline{\dim}_B \overline{O(T)}$ and the upper box dimension $\overline{\dim}_B \overline{O(T)}$ of the closure $\overline{O(T)}$ of the orbit $O(T)$ are not equal.
\item Construct trees using a set of rules in the spirit of Conway's Game of Life, and investigate the properties of the corresponding dynamical systems.
\end{enumerate}

\medskip
The rest of the paper is organized as follows. In Section \ref{gmm-section} we recall the necessary background on the space of pointed trees and the pseudogroup action on this space. In Section \ref{hausdorff} we recall the definitions of the Hausdorff and box dimensions for metric spaces. We prove Lemma \ref{lemma-holderembedding} in Section \ref{sec-holderembedding}, 
Proposition \ref{prop-subsets} in Sections \ref{non-zero} and \ref{sec-subsets}, and Proposition \ref{prop-main3}  in Section \ref{sec-application}.

\section{Pseudogroup action on the space of pointed trees} \label{gmm-section}

Let $F_n$ be a free group on $n$ generators, acting on itself on the right, and let $F_n^{1+}$ denote a non-symmetric set of generators of $F_n$. That is, if $h \in F_n^{1+}$ then $h^{-1} \notin F_n^{1+}$. Let $\cG_n$ be the Cayley graph of $F_n$. More precisely, the set of vertices is $V(\cG_n) = F_n$, and to each pair $g_1,g_2 \in V(\cG_n)$ such that $g_1 h = g_2$ for some $h \in F_n^{1+}$, we associate an edge $w_h$, labelled by $h \in  F_n^{1+}$, with $s(w_h) = g_1$ and $t(w_h) = g_2$. Thus $\cG_n$ is an oriented graph labeled by the set $ F_n^{1+}$, or, in other words, a tree where each vertex has valence $2n$. Parametrize the edges of $\cG_n$ so that each edge has length $1$, and denote by $d$ the metric on $\cG_n$ such that $d(v,w)$ is the length of the shortest path between $v$ and $w \in \cG_n$.

A subgraph $T \subset \cG_n$ is called an \emph{infinite tree} if it is non-compact and connected. Since $\cG_n$ is a tree, a connected graph $T$ is also simply connected. Denote by $X_n$ the set of all \emph{infinite pointed trees} in $\cG_n$. Such a tree is an infinite subgraph of $\cG_n$ which contains the identity $e \in F_n$. The pair $(T,e) \in X$ with the metric $d$ induced from $\cG_n$ is a pointed metric space.

\begin{defn}\label{defn-isometry}
Two (possibly finite) pointed graphs $(T,g_1)$ and $(T',g_2)$ are $F_n^{1+}$-\emph{ isomorphic} if and only if there exists an isometry $(T,g_1) \to (T',g_2)$ which maps $g_1$ onto $g_2$, such that the edge $w_h$ marked by $h$ in $T$ is mapped onto the edge $w_h'$ marked by $h$ in $T'$, with $s(w_h)$ and $t(w_h)$ mapped onto $s(w_h')$ and $t(w_h')$ respectively.
\end{defn}

Definition \ref{defn-isometry} states that two pointed graphs $(T,g_1)$ and $(T',g_2)$  are $F_n^{1+}$-isomorphic if and only if they have the same patterns of labelled edges relative to the basepoints $g_1$ and $g_2$. We remark that if $g_1 = g_2$, then an $F_n^{1+}$-isomorphism is the identity map of graphs. Indeed, in this case it preserves the distinguished vertex and the patterns of labelled edges. As a consequence it must preserve the vertex sets in $T$ and $T'$. If $g_1 \ne g_2$, then the vertex sets of $T$ and $T'$ differ even if the patterns of labelled edges are the same. In the rest of the paper we mostly consider the case $g_1 = g_2 = e$ for finite subgraphs of infinite non-isomorphic graphs. Since we treat these infinite graphs as distinct metric spaces, the use of the term $F_n^{1+}$-isomorphism for the maps of their subgraphs seems more appropriate than the use of the term the identity map.

The space $X_n$ of pointed trees is given a metric as follows. Let $B_T(e,r)$ denote a closed ball in $T$ about $e$ of radius $r$. We call such a ball in the tree $T$ a `pattern of radius $r$'.
 
\begin{defn}\label{defn-ballmetric}\cite{Ghys1999,ALM2009a,Lozano2009}
Let $X_n$ be the set of all infinite pointed trees in a locally compact Cayley graph $\cG_n$. Let $T,T' \in X_n$, and define the distance between $T$ and $T'$ by
  \begin{align} \label{eq-metric1} d_n (T,T') & = e^{-r(T,T')}, &   r(T,T') & = \max \{ r \in \mR \cup \{0\} ~ |~ \exists ~  F_n^{1+}{\rm -isomorphism} ~ B_T(e,r) \to  B_{T'}(e,r) \}. \end{align}
\end{defn}
Thus, given $T$ and $T'$ in $X$, we compare them on patterns around the basepoint $e$, and the distance between $T$ and $T'$ depends on the maximal radius of the pattern on which the corresponding subgraphs of $T$ and $T'$ are $ F_n^{1+}$-isomorphic. It is straightforward to verify that $d_n$ satisfies all axioms of a metric. In fact, \eqref{eq-metric1} satisfies the stronger form of the triangle inequality, namely
  $$d_n(T_1,T_3) = \max\{d_n(T_1,T_2),d_n(T_2,T_3)\} \ \textrm{ for any }\ T_1,T_2,T_3 \in X_n,$$
which makes $d_n$ an \emph{ultrametric}. We call $d_n$ the \emph{ball metric} on $X_n$.

 \begin{figure}[!htbp] 
\centering
    \includegraphics[width=0.7\textwidth]{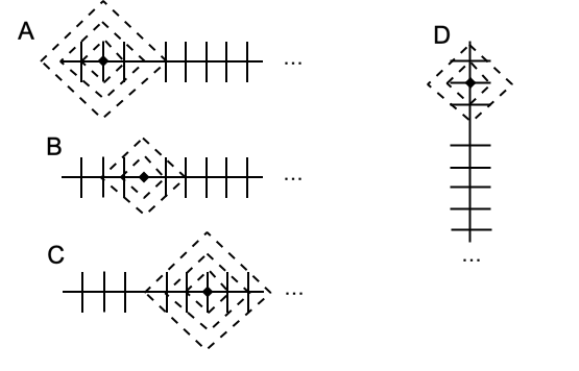}
\caption{Comparing edge patterns in metric balls of integer radius in pointed trees in $\cG_2$. For simplicity we picture the graphs which admit embeddings into the plane with edges of equal length.}    
\label{fig:metric}
\end{figure}

Definition \ref{defn-ballmetric} is illustrated in Figure \ref{fig:metric} for the case $n=2$. Denote by $a$ and $b$ the generators of a free group $F_2$. In Figure \ref{fig:metric}, the pointed trees $A$, $B$ and $C$ are infinite subgraphs of $\cG_2$ which continue infinitely far to the left and so their vertex sets contain all powers $a^m$, $m = 1,2, \ldots$ The tree $D$ continues infinitely far downward, and so it contains all powers $b^{-m}$, $m = 1,2,\ldots$ The pointed trees $A$, $B$ and $C$ have the same geometry as graphs with labelled edges, but they differ by the position of the distinguished vertex $e$. Their vertex sets in $F_2$ have a large intersection but they do not coincide. The graph $D$ has the same geometry as $A$, $B$ and $C$ as an unlabelled graph, but different geometry than $A$, $B$ and $C$ as a graph with labelled edges. For simplicity, we chose for this illustration the subtrees of $\cG_2$ which admit embeddings into the plane with edges of equal length (an embedding into the plane with this property need not exist for an arbitrary subgraph of $\cG_2$).

Dashed lines in Figure \ref{fig:metric} denote the boundaries of the closed balls of integer radius around the vertex $e$. Considering the edge patterns inside these metric balls, we see that there is no $F_2^{1+}$-isomorphism between $B_B(e,1)$ and any of the graphs $B_A(e,1)$, $B_C(e,1)$ and $B_D(e,1)$. Therefore,
  $$d_2(B,A) = d_2(B,C) = d_2(B,D) = e^0 = 1.$$
Next, there is a $F_2^{1+}$-isomorphism between $B_A(e,2)$ and $B_C(e,2)$, since these two balls contain the same patterns of labelled edges. The edge patterns in $B_A(e,3)$ and $B_C(e,3)$ differ, so there is no $F_2^{1+}$-isomorphism between $B_A(e,3)$ and $B_C(e,3)$. We conclude that
  $$d_2(A,C) = e^{-2}.$$
Finally, $B_D(e,1)$ contains the same pattern of edges as $B_A(e,1)$ and $B_C(e,1)$, while the edge pattern in $B_D(e,2)$ differs from $B_A(e,2)$ and $B_C(e,2)$. Therefore,
  $$d_2(A,D) = d_2(C,D) = e^{-1}.$$

For an integer $r>0$ denote by $D_{n}(T,e^{-r})$ an open subset of diameter $e^{-r}$ about $(T,e)$, that is,
  \begin{align*}D_{n}(T,e^{-r}) & = \left\{ (T',e) ~ | ~ d_{n}(T,T') < e^{-r} \right\}. \end{align*}

\begin{remark}\label{rem-integer}
{\rm
In most computations which use the metric $d_n$ in Definition \ref{defn-ballmetric}, one can restrict to integer values of $r$ for the following reason. Let $r>0$ and let $m = \lceil r \rceil$ be the smallest integer which is greater or equal to $r$. Consider balls of radius $r$ in trees $T$ and $T'$. The boundaries of $B_T(e,r)$ and $B_{T'}(e,r)$ intersect the edges in $T$ and $T'$ in their interior, and there is a $ F_n^{1+}$-isomorphism $B_T(e,r) \to  B_{T'}(e,r)$ if and only if there is a $ F_n^{1+}$-isomorphism $B_T(e,m) \to  B_{T'}(e,m)$. Using this property, we restrict to integer radii when computing distances between the trees in $X_n$.

Another consequence is that if $r,r'>0$ are such that $\lceil r \rceil = \lceil r' \rceil$, then the following sets are equal for any tree $T \in X$
  \begin{align}D_n(T,e^{-r}) = \{T' \mid d_n(T',T) < e^{-r}\}  = \{T' \mid d_n(T',T) < e^{-r'}\} = D_n(T,e^{-r'}) , \end{align}
and, for $m = \lceil r \rceil$,
   \begin{align}\label{eq-clopen} \{T' \mid d_n(T',T) < e^{-r}\}  = \{T' \mid d_n(T',T) \leq e^{-r}\}  = \{T' \mid d_n(T',T) \leq e^{-m}\}  , \end{align}
in particular, an open ball $D_n(T,e^{-m})$ is also closed.  We also note that the diameter of $D_n(T,e^{-m})$ is equal to it's radius $e^{-m}$, since the distance between any two trees in $D_n(T,e^{-m})$ is $e^{-m}$.
}
\end{remark}

\begin{remark}\label{rem-nofinitegraphs}
{\rm
In Definition \ref{defn-ballmetric} finite subgraphs of $\cG_n$ were excluded from consideration, since any finite subgraph $\cG_n$ is an isolated point with respect to the ball metric $d_n$, and no interesting dynamics would happen on these isolated points.
}
\end{remark}

One can show \cite{Ghys1999,Lozano2009,Blanc2001} that the metric space $(X_n,d_n)$ is compact and totally disconnected. Since we assume that the number of generators $n \geq 2$, the space $X_n$ is perfect, and so is a Cantor set.

There is a partial action of $F_n$ on the space of pointed trees $X_n$, defined as follows.

\begin{defn}\label{defn-actionany}
Let $F_n$ be a free group with non-symmetric generating set $F_n^{1+} $, and let $(X_n,d_n)$ be the corresponding set of pointed trees. Then for $g \in F_n$,
\begin{enumerate}
\item $(T,e) \cdot g$ is defined if and only if $g \in V(T)$, where $V(T) \subset F_n$ is the set of vertices of $T$.
\item $(T',e) = (T,e) \cdot g$ if and only if there is an $F_n^{1+} $-isomorphism $\alpha:(T,g) \to (T',e)$. 
\end{enumerate}
\end{defn}
In this definition $(T,e) \cdot g = (T',e)$ if and only if $T'$ has the same pattern of labeled edges relative to $e$ as the pattern of labeled edges in $T$ relative to $g$.

To the partial action of $F_n$ on $X_n$ in Definition \ref{defn-actionany} we can associate a pseudogroup of local homeomorphisms $\G_n$ as follows. 
For each $g \in F_n$ let $\ell_g = d(e,g)$, where $\ell_g$ is the distance between $e$ and $g$ in the word metric on the group $F_n$. There is a finite number of subgraphs of $\cG_n$ of radius $\ell_g$, which contain $g$, so by Property $(1)$ in Definition \ref{defn-actionany} the action of $g$ is defined on the finite union of clopen subsets
  \begin{align}\label{domgamma}D & = \bigcup \left\{D_n(T,e^{-\ell_g}) ~ |~ T \in X_n, ~ g \in V(T) \right\}. \end{align}
  
The mapping
  \begin{align*} 
  \gamma_g & : {\rm dom}(\gamma_g) = D \to X_n : (T,e) \mapsto (T,e)\cdot g 
  \end{align*}
is a homeomorphism onto its image, and a pseudogroup $\G_n$ is defined to be a collection  
  \begin{align}\label{eq-pseudogroup}\G_n& = \bigsqcup \{ ~\gamma_g ~|~ g \in F_n ~ \}\end{align}
of local homeomorphisms (see \cite{CandelConlon2000} for more details about pseudogroups). 
The subset $\G^1_n = \{~\gamma_g \in \G_n ~|~ g \in F^{1+}_n ~\}$ is a generating set of $\G_n$. We remark that a pointed tree $(T,e)$ is in the domain of $\gamma_g$ if and only if $T$ contains the vertex $g$ if and only if $T$ contains the unique finite path of edges joining the vertices $e$ and $g$ in the Cayley graph $\cG_n$ of the free group $F_n$.

\begin{remark}  \label{rem-lipschitz}
{\rm
An easy computation shows that the generators in $\G_n^1$ are Lipschitz maps. Indeed, if $T,T' \in {\rm dom}(\gamma_g)$, $g \in F^{1+}_n$, and $d_n(T,T') = e^{-m}$, then 
   \begin{align*} d_n(\gamma_g(T),\gamma_g(T')) \leq e^{-(m-1)}, \end{align*}
and so $\gamma_g$ is Lipschitz with $C = e$.
}
\end{remark}

\begin{defn}\label{defn-eqrel}
Let $(X_n,d_n)$ be the  space of pointed trees, and $\G_n$ be the pseudogroup on $X_n$. Then the orbit of a graph $(T,e) \in X_n$ is the set 
  $$O(T) = \{T' \in X_n | \textrm{ there exists }\gamma_g \in \G \textrm{ such that }\gamma_g(T,e) = (T',e)\} .$$
\end{defn}

In Definition \ref{defn-eqrel}, the tree $T'$ has the same pattern of edges relative to $e$ as $T$ has relative to the vertex $g$.

\section{The box and the Hausdorff dimensions}\label{hausdorff}

In this section, we define the  Hausdorff dimension, and the lower and upper box dimensions, for the metric space of pointed graphs $(X_n,d_n)$, associated to the free group $F_n$, for $n \geq 2$.  See \cite{Falconer1985,Pesin1997} for further details and the properties of these dimensions.

Let $\cU$ be the collection of all open subsets in the space of pointed graphs $(X_n,d_n)$, $n \geq 2$. Let $Z \subset X_n$ be a subset of $X_n$, and let the numbers $\alpha>0$ and $\e>0$ be given. Since $X_n$ is compact, there always exists a finite subcollection $\cV \subset \cU$ of open sets of diameter at most $\e$, whose union contains $Z$ (since the closure of $Z$ is compact in $X_n$). Therefore, we can define a function
   \begin{align}\label{h-one} M_H(Z,\alpha,\e) = \inf_\cV \left\{ \sum_{V \in \cV} (\diam \ V)^\alpha\right\}, \end{align}
 where the infimum is taken over all finite or countable subcollections $\cV$ of open sets of diameter at most $\e$ which cover $Z$. We define
  \begin{align}\label{h-two} m_H(Z,\alpha) = \lim_{\e \to 0} M_H(Z,\alpha,\e). \end{align}  
Since $M_H(Z,\alpha, \e)$ is non-decreasing when $\e$ decreases, this limit exists or is infinite. By \cite[Proposition 1.2]{Pesin1997}, we have that for a set $Z$, there exists a critical value $\alpha_H$, called the \emph{Hausdorff dimension} of $Z$, such that 
  \begin{align}\label{h-three} \dim_H Z =\alpha_H = \inf\left\{ \alpha: m_H(Z,\alpha) = 0 \right\} = \sup \left\{ \alpha : m_H(Z,\alpha) = \infty \right\}.\end{align}

The definition of the Hausdorff dimension uses covers by arbitrary open sets with bounded diameters. In actual computations it may be easier to deal with \emph{open balls} around a point in the space. 
It is an exercise to show that, for a metric space, the Hausdorff dimension of sets does not change if one considers a collection of open balls $\cB$ in $X_n$ instead of all open sets. So in the rest of the paper, we consider the collection of \emph{open balls} instead of arbitrary open sets. By Remark \ref{rem-integer}, open balls in $(X_n,d_n)$ are also closed.

For basic properties of Hausdorff dimension we refer to \cite[Theorem 6.1]{Pesin1997}. Now let
 \begin{align}\label{b-one} R(Z,\alpha,\e) = \inf_\cV \left\{ \sum_{V \in \cV} (\diam \ V)^\alpha\right\}, \end{align}
 where the infimum is taken over all finite or countable subcollections $\cV$ of open sets \emph{of diameter precisely} $\e$ which cover $Z$. We define
  \begin{align}\label{b-two} \underline{r}(Z,\alpha) = \liminf_{\e \to 0} R(Z,\alpha,\e), && \overline{r}(Z,\alpha) = \limsup_{\e \to 0} R(Z,\alpha,\e). \end{align}  
By \cite[Proposition 1.2]{Pesin1997}, there exist numbers
  \begin{align*} \underline{\dim}_B Z = \inf \{\alpha \mid \underline{r}(Z,\alpha) = 0 \}=   \sup \{\alpha \mid \underline{r}(Z,\alpha) = \infty\} \\
  \overline{\dim}_B Z = \inf \{\alpha \mid \overline{r}(Z,\alpha) = 0 \}=   \sup \{\alpha \mid \overline{r}(Z,\alpha) = \infty\},\end{align*}
called the \emph{lower} and the \emph{upper box dimensions} of the set $Z$ respectively. Recall that we always have $\underline{\dim}_B Z \leq \overline{\dim}_B Z$. For any $\e>0$ and any $Z \subset X_n$, denote
  \begin{align}\label{Lambda-box} \Lambda(Z,\e) = \inf_\cV \{ {\rm card} \{ \cV\} \}, \end{align}
where $\cV$ runs over all finite or countable subcollections $\cV$ of open sets of diameter precisely $\e = e^{-m}$ which cover $Z$. Then by \cite[Theorem 2.2]{Pesin1997} we have
  \begin{align} \label{eq-boxcompute1} \underline{\dim}_B Z = \liminf_{\e \to 0} \frac{\log \Lambda (Z,\e)}{ \log (1/\e)} = \liminf_{m \to \infty} \frac{1}{m} \log \Lambda(Z, e^{-m}),  \\   \label{eq-boxcompute2}\overline{\dim}_B Z = \limsup_{\e \to 0} \frac{\log \Lambda (Z,\e)}{ \log (1/\e)} = \limsup_{m \to \infty} \frac{1}{m} \log \Lambda(Z, e^{-m}).
  \end{align}

\section{H\"older embeddings of shift spaces}\label{sec-holderembedding}

 In this section, we consider the Hausdorff dimensions of the images of the  embeddings of the standard shift spaces into the pseudogroup dynamical systems $(X_n,\G_n)$, for $n \geq 2$,  as constructed by Lozano Rojo and the author in \cite{LozanoLukina2012}.

Let $\cA$ be a finite alphabet, and let $\Sigma = \{\sigma: \mZ \to \cA\}$ be the space with the shift action of $\mZ$;  that is, for $k \in \mZ$,
  $$\sigma (m) \cdot k = \sigma (k+m), \textrm{ for all } m\in \mZ.$$
A basis for the topology on $\Sigma$ is formed by the cylinder sets 
   \begin{align}\label{eq-cylinder} [s_0,s_1,\ldots, s_\ell]_{t} = \{\sigma \in \Sigma \mid \sigma(t) = s_0, \ldots, \sigma(t+\ell) = s_\ell \}.  \end{align}
Recall the construction of the embedding of $\Sigma$ into $X_n$, $n \geq 2$, from \cite{LozanoLukina2012}.

Let $n \geq {\rm card}(\cA)$. Choose an injective map $\alpha: \cA \to F_n^{1+}$, where $F_n^{1+}$ is the `positive part' of the symmetric generating set $F^1_n$ of $F_n$. More precisely, if $g \in F^{1+}_n$, then $g^{-1} \notin F^{1+}_n$. Construct a map
  \begin{align}\label{eq-embedmap} \Phi_\alpha : \Sigma \to X_n \end{align} 
  inductively as follows. 
  
\begin{figure}[!htbp] 
      \noindent
\begin{minipage}{.32\textwidth}
\centering
    \includegraphics[width=0.7\textwidth]{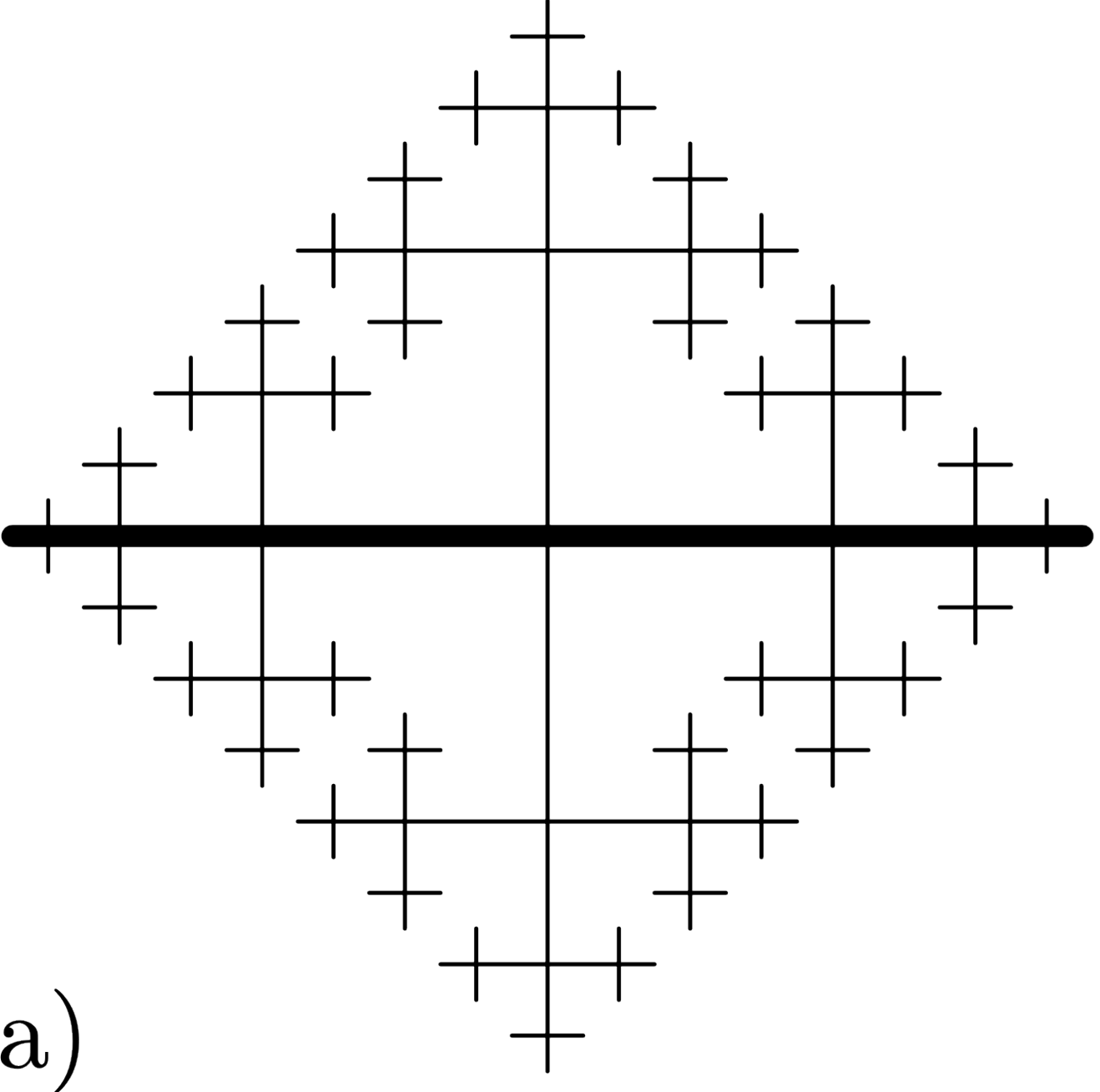}
 \end{minipage}
 \begin{minipage}{.32\textwidth}
\centering
    \includegraphics[width=0.7\textwidth]{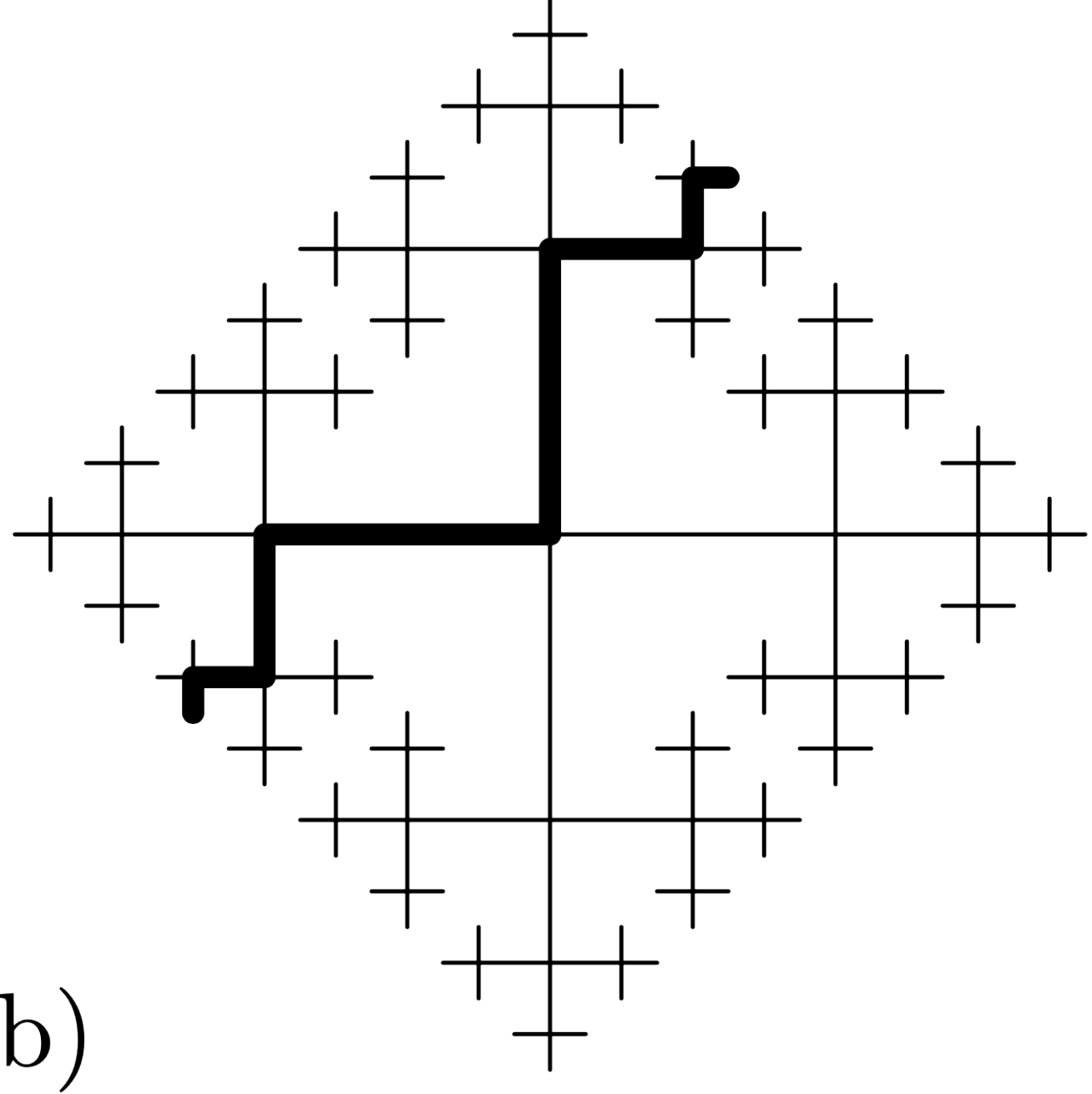}
 \end{minipage}
 \begin{minipage}{.32\textwidth}
\centering
    \includegraphics[width=0.7\textwidth]{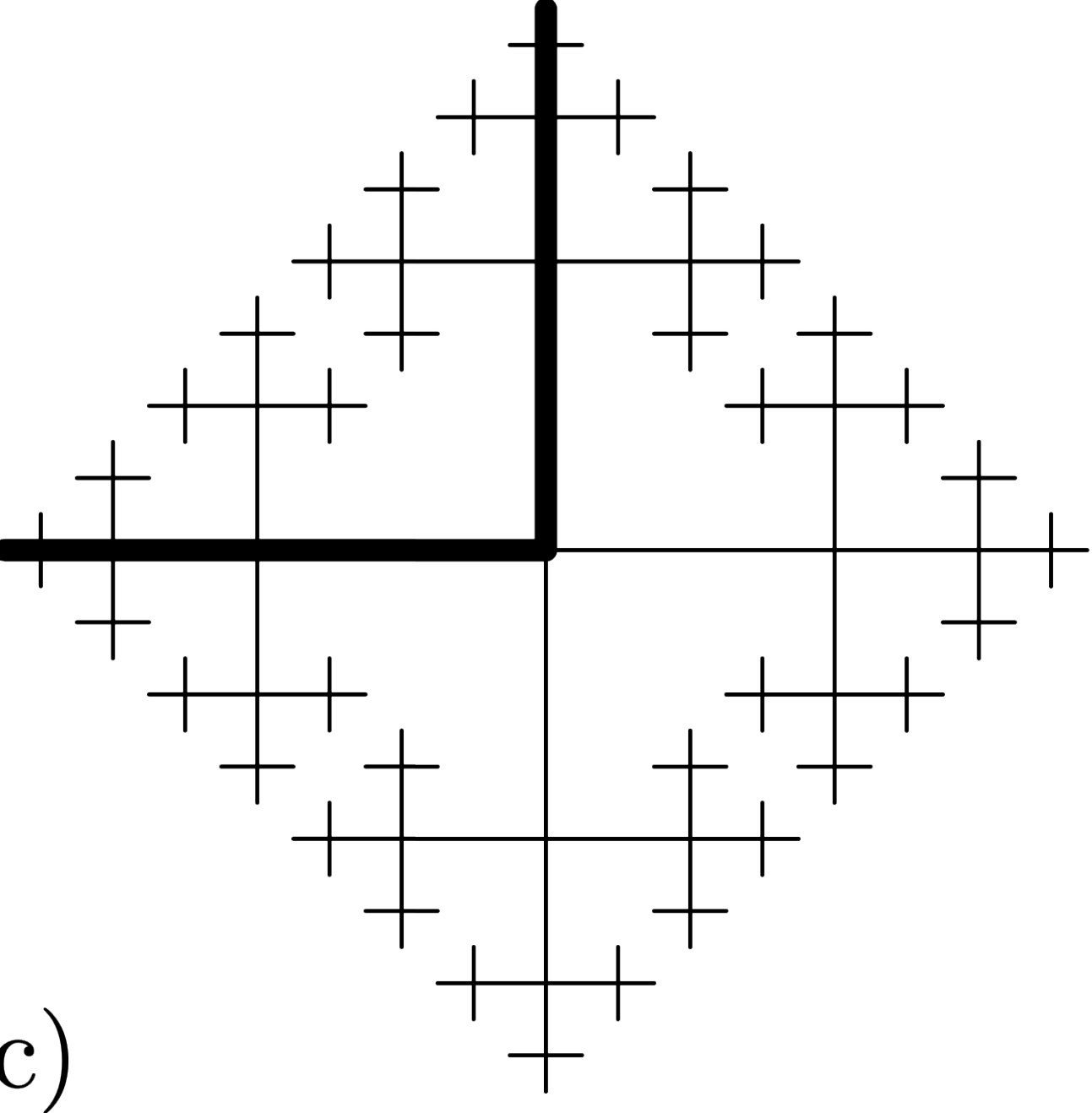}
 \end{minipage}
   \caption{ For $F_2^+ = \langle a,b \rangle$, where $a$ labels horizontal edges, and $b$ labels vertical edges, $\mathcal{A} = \{0,1\}$ and $\alpha(0) = b$, $\alpha(1) = a$, the tree $T_{\sigma}$ for specific $\sigma$ is shown as a subgraph of $\cG_2$. Here: a)   $\sigma = (...1111.1111 ...)$;  b) $\sigma = (...0101.0101 ...)$; c) $\sigma = (...1111.0000 ...)$.}
    \label{fig:shiftembedded}
  \end{figure}

  Let $\sigma \in  \Sigma$, and denote $\sigma_k = \sigma \cdot k$. Set $\mZ^0 = \{0\}$, and $\mZ^{ \leq i} = \{-i, -(i-1), \ldots, i-1, i\}$. Let  $K_0 = \{e\}\subset\cG_n$ be a graph consisting of a single vertex $e$ and no edges, and define 
  $$\lambda_0:\mZ^{\leq 0} = \mZ^0 \to V(K_0): 0 \mapsto e.$$
By induction, suppose we have a connected subgraph $K_{i-1} \subset \cG_n$ containing $e$, and a bijection $\lambda_{i-1} : \mZ^{\leq i-1} \to V(K_{i-1})$ such that for $-(i-1)<k \leq i-1$ we have
 \begin{align*} \lambda_{i-1}(k) = \lambda_{i-1}(k-1) \alpha(\sigma_{k-1}(e)) ,\end{align*}
where the juxtaposition denotes the group operation in $F_n$. 
 Then obtain a map 
   $$\lambda_i:\mZ^{\leq i} \to V(K_i)$$
 as follows. Note that $\mZ^{\leq i-1} \subset \mZ^{\leq i}$, and for any $k \in \mZ^{\leq i-1}$ set $\lambda_i(k) = \lambda_{i-1}(k)$. Then define
  $$ \lambda_i(i) = \lambda_{i-1}(i-1) \alpha(\sigma_{i-1}(e)),  \textrm{ and }  \lambda_i(-i) = \lambda_{i-1}(-(i-1)) (\alpha(\sigma_{-i}(e)))^{-1}.$$
 Then let $K_i$ be a subgraph of $\cG_n$ with the set of vertices $V(K_i) = \lambda_i(\mZ^{\leq i})$. Then define 
  \begin{align}\label{eq-themap}\Phi_\alpha(\sigma) = T_\sigma = \bigcup_{i=0}^\infty K_i. \end{align}
 It is proved in \cite{LozanoLukina2012} that the map \eqref{eq-embedmap} is an equivariant embedding. This embedding depends on the choice of the map $\alpha$. Since $\mZ$ has one generator, for every tree $T_\sigma \in \Phi_\alpha(\Sigma)$ the vertices in $V(T_\sigma)$ have valence two, that is, every vertex is adjacent to precisely $2$ edges, one outgoing, and one incoming, see Figure \ref{fig:shiftembedded}. The labels of the edges are the generators in the image of $\alpha$ in $F_n^{1+}$. 
  
  A metric on $\Sigma$ generating the topology, can be given, for example, by (see for instance \cite{Kitchens})
  \begin{align*} d(\sigma_1,\sigma_2) = {2^{-k}}, \textrm{ where } k = \max \{ m \mid \sigma_1(i) = \sigma_2(i) \textrm{ for } |i| \leq m\}. \end{align*}
 It is well-known that the Hausdorff and the box dimensions of a shift space are positive and equal.
   
 We now restate and prove Lemma \ref{lemma-holderembedding}.
 
\begin{lemma}\label{shift}
Let  $\cA$ be a finite alphabet, and $\Sigma = \{\sigma: \mZ \to \cA\}$ be the shift space. Suppose there is an equivariant embedding $\Phi: \Sigma \to X_n$ as in \cite{LozanoLukina2012}. Then $\Phi$ is H\"older continuous, and, moreover,
  \begin{align}\label{eq-shiftdim}0<\dim_H \Phi(\Sigma) = \underline{\dim}_B \Phi(\Sigma) =  \overline{\dim}_B \ \Phi(\Sigma) < \infty.\end{align}
  \end{lemma}

\proof
The embedding $\Phi:(\Sigma,d) \to (X_n,d_n)$ is H\"older continuous if there exist constants $C>0$ and $0 < \gamma < 1$ such that for any $\sigma_1,\sigma_2 \in \Sigma$ we have 
  \begin{align*}   d_n(\Phi(\sigma_1), \Phi(\sigma_2)) \leq C d(\sigma_1,\sigma_2)^{\gamma}. \end{align*}
Notice that $d(\sigma_1,\sigma_2) = {2^{-k}}$ if and only if $d_n(\Phi(\sigma_1),\Phi(\sigma_2)) = e^{-k}$. 
Then $\Phi$ is H\"older continuous with $C = 1$ and $\gamma = 1/\ln 2$ and, in fact, a stronger statement than just H\"older continuity is true:
  $$d_n(\Phi(\sigma_1), \Phi(\sigma_2)) = d(\sigma_1,\sigma_2)^{1/\ln 2}.$$
We note that if $V \subset \Sigma$ is open, then the properties of the ultrametrics $d$ and $d_n$ imply that $\Phi(V) = \Phi(\Sigma) \cap \widehat{V}$ for a clopen set $\widehat V$ such that 
  $$\diam(\widehat V) = \diam(\Phi(V)) = (\diam(V))^{1/ \ln 2}.$$
Denote by $\widehat \cV$ an open cover of $\Phi(\Sigma)$ in $X_n$ by sets of diameter $\e$, and by $\cV$ the corresponding cover of $\Sigma$ by sets of diameter $\e^{\ln 2}$. Using the above equality of diameters of sets in formula \eqref{Lambda-box} we obtain
 \begin{align*} \Lambda(\Phi(\Sigma),\e) = \inf_{\widehat \cV} \{ {\rm card} \{ \widehat \cV\} \} = \inf_{\cV} \{{\rm card} \{ \cV\}\} = \Lambda(\Sigma,\e^{\ln 2}).\end{align*}
Then for the lower box dimension \eqref{eq-boxcompute1} gives
$$\underline{\dim}_B \Phi(\Sigma) = \liminf_{\e \to 0} \frac{\log \Lambda (\Phi(\Sigma),\e)}{ \log (1/\e)} =  \liminf_{\e \to 0} \frac{\log \Lambda (\Sigma,\e^{\ln 2})}{ \log (1/\e^{\ln 2})^{1/\ln 2}} = \ln 2 \, \underline{\dim}_B \Sigma.$$
By a similar computation we obtain 
  $$\overline{\dim}_B \Phi(\Sigma) = \ln 2 \ \overline{\dim}_B \Sigma  = \ln 2\ \underline{\dim}_B \Sigma = \underline{\dim}_B \Phi(\Sigma).$$
For the Hausdorff dimension, since $\Phi$ is H\"older, we obtain
   $$\dim_H \Phi(\Sigma) \leq \ln 2 \dim_H \Sigma < \infty.$$
Next, the inverse $\Phi^{-1}: \Phi(\Sigma) \to \Sigma$ from the image of the embedding to $\Sigma$ is H\"older continuous with $C = 1$ and $\gamma = \ln 2$, and we obtain
    $$0< \dim_H \Sigma \leq 1/\ln 2 \dim_H \Phi( \Sigma),$$
 which shows that $\dim_H \Phi( \Sigma) = \ln 2 \dim_H \Sigma =  \ln 2 \ \overline{\dim}_B \Sigma = \overline{\dim}_B \, \Phi(\Sigma) $.   This proves \eqref{eq-shiftdim}.
\endproof

\section{Hausdorff dimension of the space of pointed trees} \label{non-zero}

 In this section, we show that   $(X_n, d_n)$  has infinite Hausdorff dimension. We give the proof for $n=2$ in Lemma \ref{thm-matchbox}. The proof for $n>2$ is obtained as a consequence of Lemma \ref{thm-matchbox} in Corollary \ref{cor-matchbox}.

Let us first fix some terminology. In the following, we use two types of metric spaces: graphs $T \subset \cG_2$ with length metric denoted by $d$, and the space $X_2$ with the ultrametric $d_2$.  

Let $B_T(e,m)$ denote the \emph{closed} ball of radius $m$ in the graph $T$ about the basepoint $e$, with respect to metric $d$ on $T$. We say that $B_T(e,m)$ is a \emph{pattern} defined by $T$, as this is just a subgraph of $T$. The subscript $T$ in the notation indicates that $B_T(e,m)$ is a pattern determined by $T$.

We reserve the word `ball' for clopen balls in $X_2$, where $D_2(T,e^{-m})$ denotes the ball of radius $e^{-m}$ about $T$ in the space $X_2$. Such a ball contains all pointed trees $T'$, which contain the same pattern of radius $m$ as $T$ relative to the basepoints $e$. Since $d_2$ is an ultrametric, the diameter of $D_2(T,e^{-m})$ is also $e^{-m}$.

\begin{prop}\label{thm-matchbox}
Let $F_2$ be a free group on $2$ generators, and let $(X_2,d_2)$ be the corresponding space of pointed graphs with the ball metric. Then  $\dim_H(X_2) = \infty$.
\end{prop}

\proof Let $\cU$ be a finite cover of $X_2$ by clopen balls of diameter at most $e^{-m}$. Since $d_2$ is an ultrametric, two clopen sets $U \cap W$ have non-empty intersection if and only if one of them is contained in another one, therefore, we can assume that $\cU$ is a partition of $X_2$.

For a set $Z \subset X_2$, such that  $T \in Z$ if and only if each vertex $v \in V(T)$ has precisely $3$ adjacent edges, we are going to show that $\dim_H Z = \infty$. Then by \cite[Theorem 6.1]{Pesin1997}, since $Z \subset X_2$, then $\dim_H X_2 = \infty$. 

Given a pattern $P$ of radius $i$, such that every interior vertex in $P$ has precisely $3$ adjacent edges, and a tree $T_P \in Z$ containing $P$, we estimate the number of clopen balls of diameter $e^{-(i+1)}$ which are needed to cover a clopen ball $D_2(T_P,e^{-i})$ about $T_P$.

\begin{figure}[!htbp]
  \noindent
\begin{minipage}{.4\textwidth}
\centering
\includegraphics[width=0.6\textwidth]{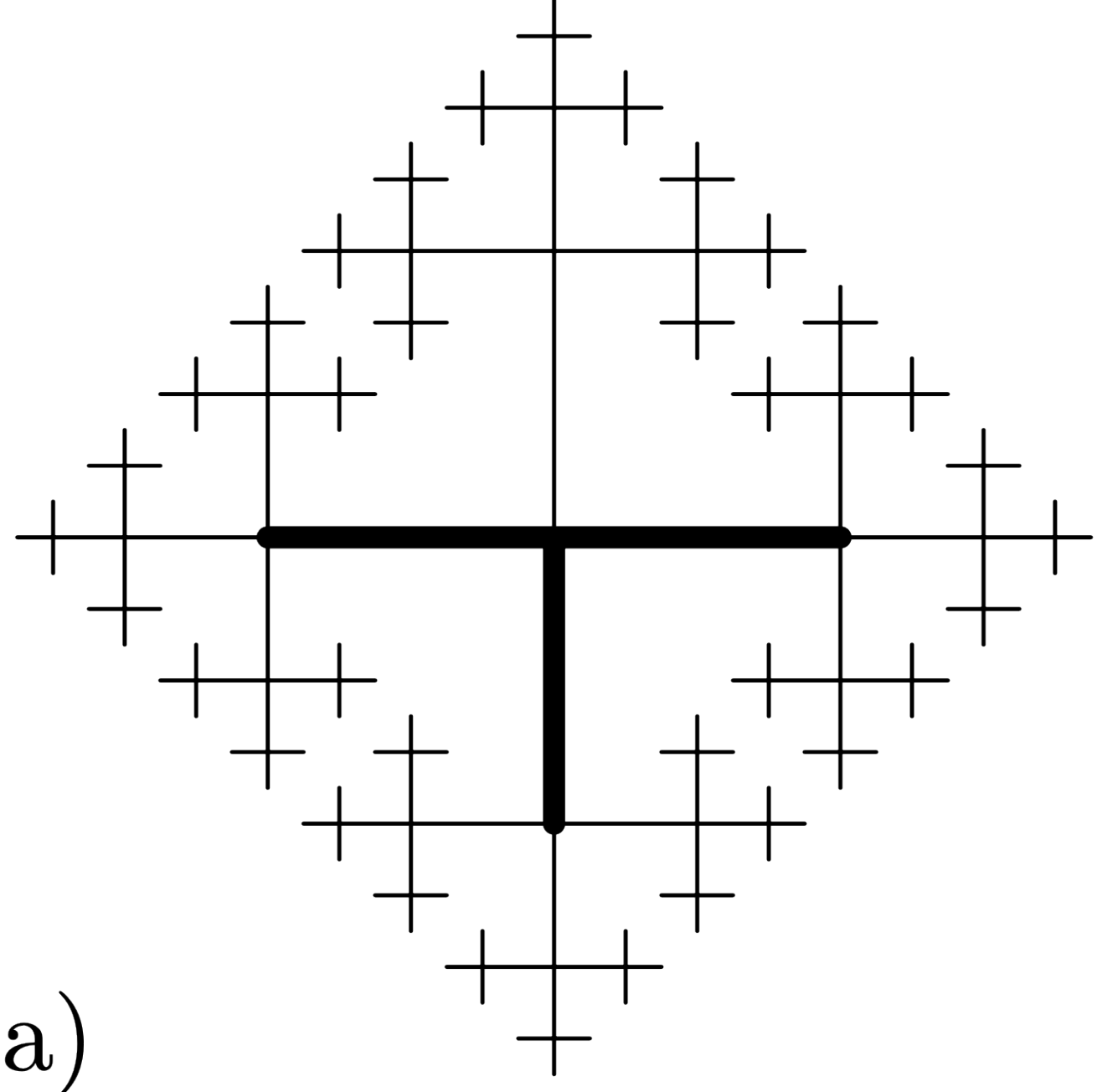}
  \end{minipage}%
  \begin{minipage}{.4\textwidth}
\centering
\includegraphics[width=0.6\textwidth]{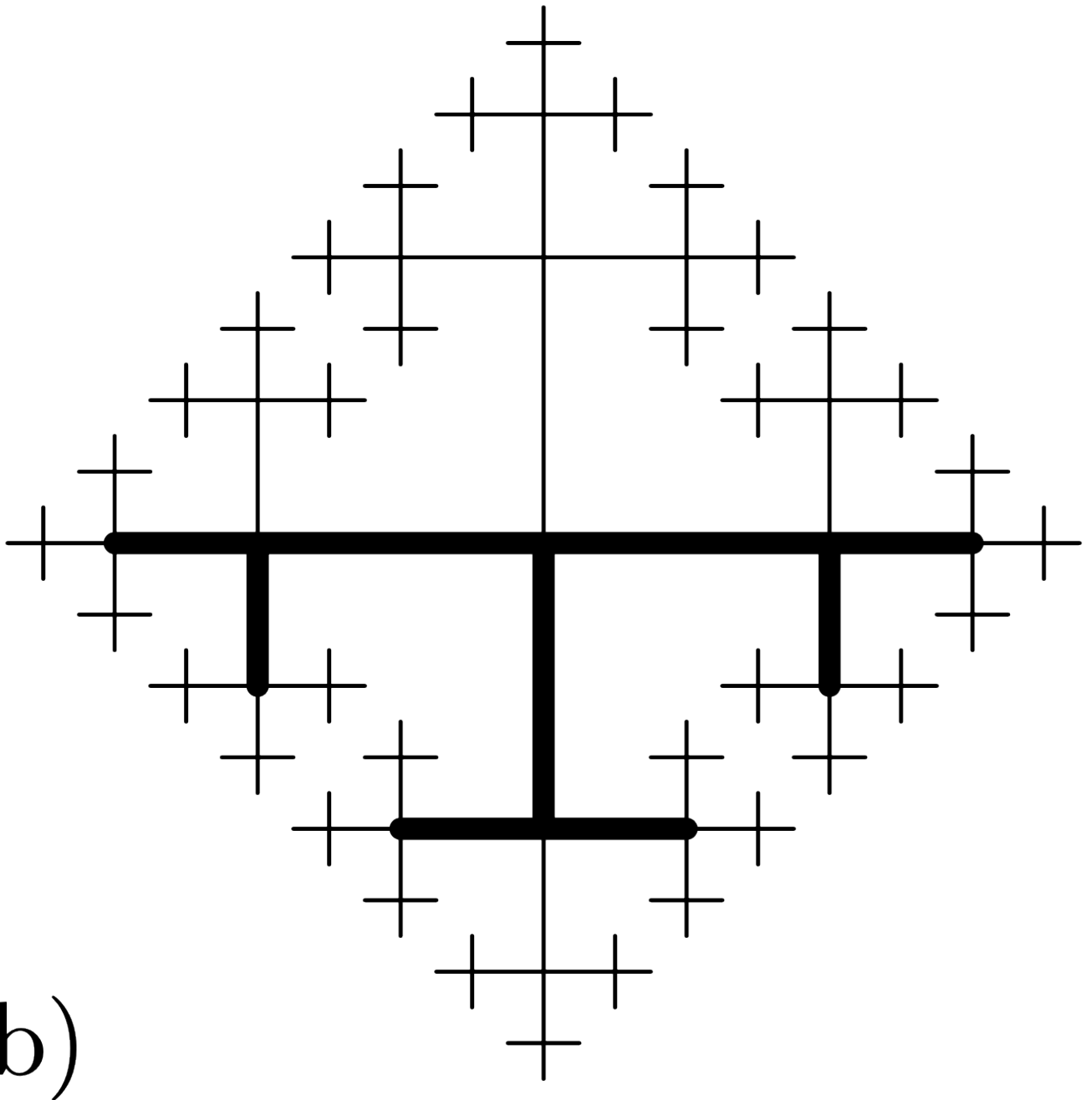}
  \end{minipage}%

    \caption{ a) A pattern of radius $1$ that has three vertices in it's $1$-boundary. b) A pattern of radius $2$ that has $3\cdot 2$ vertices in it's $2$-boundary. The patterns are shown as subgraphs of $\cG_2$.}
    \label{fig:boundaryvertices}
  \end{figure}

Denote by $\partial_i P$ the `$i$-boundary of a pattern $P$', which consists of the vertices in $P$ which are at distance precisely $i$ from the basepoint $e$, relative to the metric $d$ on $T_P$, induced from $\cG_2$. For example, the $1$-boundary of a pattern of radius $1$ inside $P$ contains $3$ vertices, the $2$-boundary of a pattern of radius $2$ inside $P$ contains $3*2$ vertices, and, inductively, the $i$-boundary of the pattern $P$ contains $K = 3 \cdot 2^{i-1}$ vertices, see Figure \ref{fig:boundaryvertices}. 

To determine how many clopen balls in $Z$ of diameter $e^{-(i+1)}$ are contained in $D_2(T_P,e^{-i})$, we consider all possible patterns $Q$ of radius $(i+1)$ which coincide with $P$ on a subset of radius $i$ around the identity vertex, and which have precisely $3$ edges attached to any interior vertex. Such a $Q$ is obtained by attaching $2$ edges to every vertex in $\partial_i P$. Since every vertex in $\cG_2$ has $4$ adjacent edges, there are $3$ distinct ways to do that for each vertex in $\partial_i P$. Then
   \begin{align}\label{partition-size} S_{P,i+1} & = 3^{3 \cdot 2^{i-1}} , \end{align}
which is the number of clopen sets of diameter $e^{-(i+1)}$ required to cover a clopen set of diameter $e^{-i}$.  We show that, since $S_{P,i+1}$ grows superexponentially with $i$, $Z$ cannot have finite Hausdorff dimension.

We note that $Z$ is a closed, and, therefore compact subset of $X_2$. Indeed, if $T$ is an accumulation point of $Z$, then any vertex in $T$ has precisely three adjacent edges, and so $T \in Z$. Thus we can restrict to finite covers of $Z$. For a finite cover of $Z$ by clopen sets of diameter at most $e^{-m}$  there is only a finite number of  set diameters which can occur in this cover. Let $ e^{-(m+s)}$ be the minimal diameter, then we can write
  \begin{align}\label{eq-cover11} \sum_{U \in \cU} (\diam ~ U)^\alpha = N_{m} e^{-m\alpha} + N_{m+1} e^{-(m+1)\alpha}  + \cdots +  N_{m+s} e^{-(m+s)\alpha} ,\end{align}
where $N_{m_i} \geq 0$ are integers. If $\dim_H Z < \infty$ is finite, for any $\alpha>\dim_H Z$ we have
  \begin{align}\label{conv-eq} \lim_{m \to \infty} \inf_{\cU_m} \sum_{U \in \cU_m} (\diam (U))^{\alpha} = 0, \end{align}
where $\cU_m$ runs through all finite partitions of $Z$ by clopen balls of diameter at most $e^{-m}$. The convergence of the limit \eqref{conv-eq} to zero implies that there exists a sequence of real numbers $0<\varepsilon_m <1$ such that $\ds \lim_{m \to \infty}\varepsilon_m  = 0$ and for any $m>0$,
   $$ \inf_{\cU_m} \sum_{U \in \cU_m} (\diam (U))^{\alpha} < \varepsilon_m. $$
Then for any $m$ there exists a finite partition $\widetilde{\cU}_m$ of $Z$ such that, for some $0<\epsilon< \varepsilon_m$ we have
  \begin{align}\label{eq-ineqcovers} \inf_{\cU_m} \sum_{U \in \cU_m} (\diam (U))^{\alpha} \leq  \sum_{U \in \widetilde{\cU}_m} (\diam (U))^{\alpha} < \inf_{\cU_m} \sum_{U \in \cU_m} (\diam (U))^{\alpha} +\epsilon < \varepsilon_m. \end{align}

The partition $\widetilde{\cU}_m$ satisfies the equation \eqref{eq-cover11}, and combining that with \eqref{eq-ineqcovers} we obtain 
   \begin{align} \label{eq-boundcover}N_{m+i} e^{-(m+i) \alpha} < \varepsilon_{m} <1, \end{align}
for each $1 \leq i \leq s$. This yields 
   \begin{align}\label{Nr-est}N_{m+i} < e^{(m+i) \alpha}, ~1 \leq i \leq s .\end{align} 
That is, the number of clopen sets $N_{m+i}$ of diameter $e^{-(m+i)}$ in $\widetilde{\cU}_m$ must grow at most exponentially with $i$. However, by  \eqref{partition-size} the number of clopen sets of diameter $e^{-(m+i+1)}$ required to cover a clopen set of diameter $e^{-(m+i)}$ grows superexponentially. It follows that the cover with property \eqref{conv-eq} cannot exist, and so $\dim_H Z = \infty$.
\endproof

Since $X_2$ can be isometrically embedded into $X_n$ for any $n \geq 2$, and isometries preserve the Hausdorff dimension of sets \cite{Pesin1997}, the following corollary is straightforward.

\begin{cor}\label{cor-matchbox}
For $n \geq 2$, the Hausdorff dimension of $X_n$ is infinite, i.e. $\dim_H(X_n) = \infty$.
\end{cor}

\section{Subsets of the space of pointed graphs}\label{sec-subsets}

In this section, we prove statements $2)$ and $3)$ of the Proposition \ref{prop-subsets}. We start with statement $2)$.

Recall that for a finite pattern $P$ of radius $\ell$ we denote by $\partial_\ell P$ the set of vertices in $P$ which are at distance $\ell$ from the identity $e$ in $P$, with respect to the metric $d$ on $P$, induced from the length metric on $\cG_n$. We call this set the `$\ell$-boundary' of $P$.

\begin{prop}\label{non-box}
The set of closed invariant transitive subsets with zero Hausdorff dimension and positive box dimension is dense in $(X_n,d_n)$.
\end{prop}

\proof Let $T_P$ be a tree in $X_n$, and let $P = B_T(e,\ell)$ be a pattern of radius $\ell$ in $T_P$. Since $T_P$ is connected and non-compact,  the $\ell$-boundary of $P$ is non-empty and contains a vertex $v \in \partial_\ell P$.

\begin{figure}[!htbp]
  \noindent
\begin{minipage}{.4\textwidth}
\centering
\includegraphics[width=0.6\textwidth]{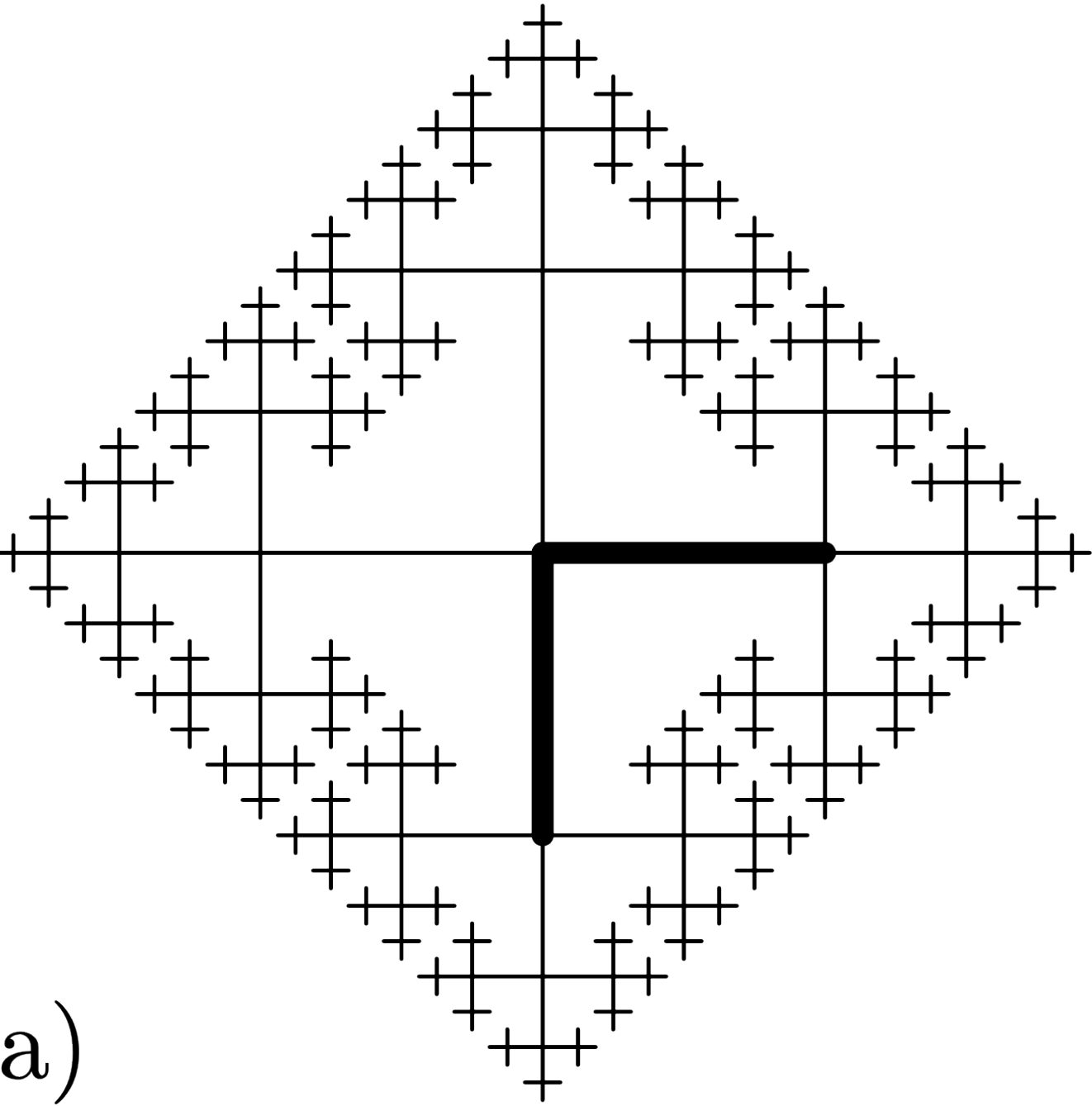}
\end{minipage}%
\begin{minipage}{.4\textwidth}
\centering
\includegraphics[width=0.6\textwidth]{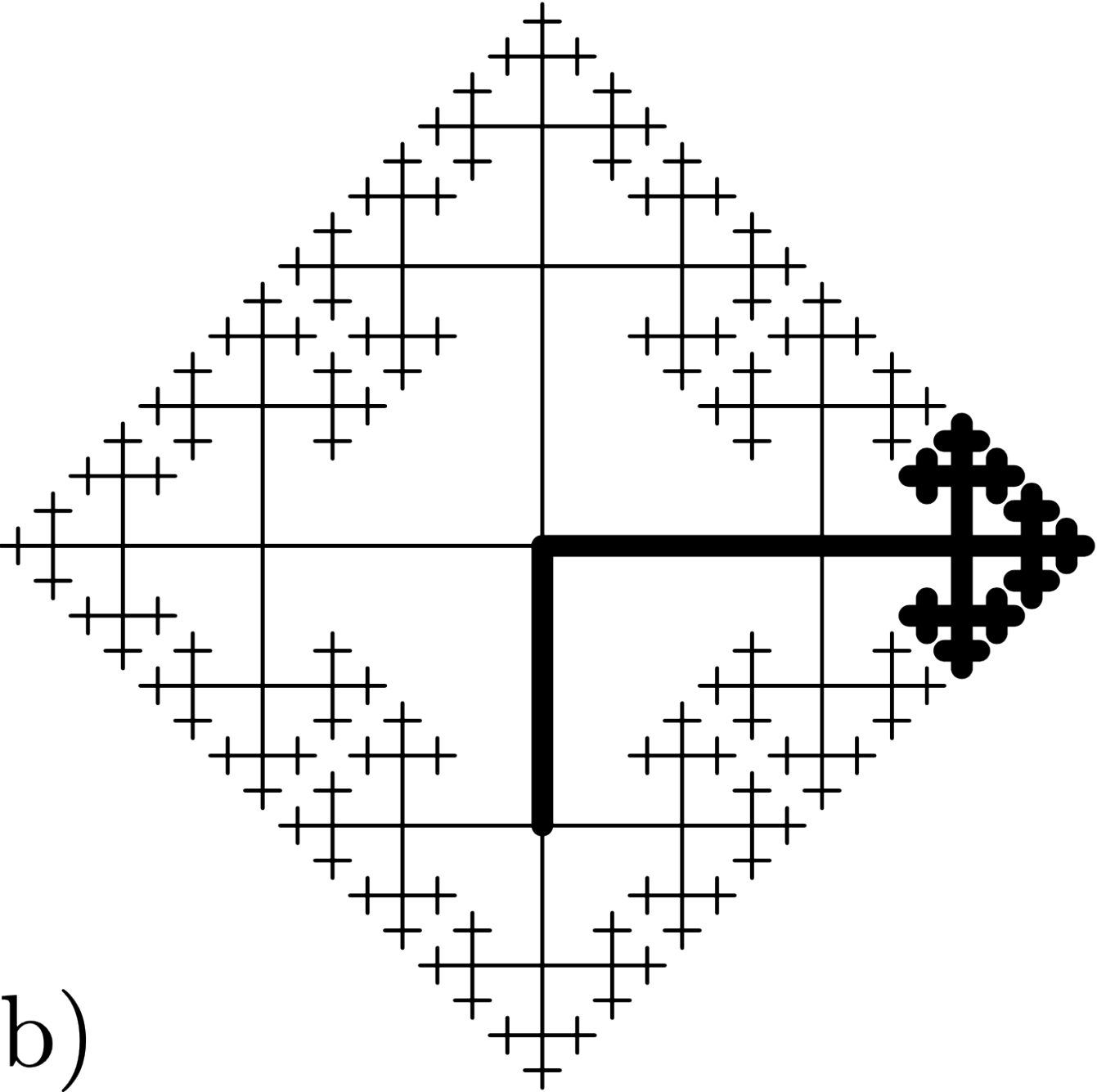}
\end{minipage}%

    \caption{ a) An example of a pattern $P$. The pattern $P$ has radius $1$, and has $2$ vertices in its $1$-boundary.  Here $h = a$, the generator marking the horizontal edges, $n=2$, the vertex $v = a$. b) The pattern $P$ with a `branch' of $\cG_2$ attached at $v=a$. The patterns are shown as subgraphs of $\cG_2$.}
    \label{fig:tail}
  \end{figure}

Construct another tree $T \supset P$ as follows. Let $e_h$ be an edge in the interior of $P$ with endpoint $v$. Set $k = 1$ if $v=t(e_h)$, and set $k = -1$ if $v = s(e_h)$. The set of vertices of $T$ is defined by
 $$V(T) =  V(P) \cup\{ v h^k mg ~| ~ m \ne h^k, ~ mg \in \cG_n {\textrm{~ is a reduced word}} ~\}.$$ 
We can think of $T$ as the union of $P$ and of an infinite `branch' of the graph $\cG_n$, starting with $h^k$, attached to $v$ instead of $e$, see Figure \ref{fig:tail} for an example. The branch of $\cG_n$ is infinite, and every vertex in the branch except $v$ is adjacent to $2n$ edges. The vertex $v$ is attached to only $2$ edges. Therefore, for any $g \in V(T)$ we have $(T,e) \cdot g \ne (T,e)$, that is, $(T,e)$ is not fixed by the action of any $g \in \G_n$.

Consider the closure $M = \overline{O(T)}$ of the orbit of $T$ under the action of the holonomy pseudogroup $\G_n$. We are going to show that $M$ contains just two orbits, and so it is countable.

Let $T' \in M$. Then for every $k>0$ there exists $v_k \in V(T)$ such that there is an isomorphism $B_{T'}(e,k) \to B_T(v_k,k)$. Let $k = 2 \ell +1$, where $\ell$ is the radius of the pattern $P$. If $B_T(v_{2 \ell +1},2\ell+1)$ does not contain the vertex $v$, then $B_T(v_{2 \ell +1},2\ell+1)$ is isomorphic to $B_{\cG_n}(e,2\ell+1)$. Then either $B_T(v_k,k)$ is isomorphic to $B_{\cG_n}(e,k)$ for all $k> 2 \ell+1$, in which case $T' = \cG_n$, or $B_T(v_k,k)$ contains $v$ for some $k > 2\ell+1$. If $v_k \in P$, then $B_T(v_k,k)$ contains $P$ and intersects the branch attached to $v$. Then $v_m = v_k \in P$ for all $m \geq 0$, and so $T' \in O(T)$. If $v_k \notin P$, then $v_k$ is in the branch of $\cG_n$, attached to $v$. Let $\lambda = d(v_k,v)$. Since all vertices in the branch of $\cG_n$ attached to $v$ are adjacent to $2n$ edges, $v$ is adjacent to $2$ edges, and all patterns $B_T(v_m,m)$ must match on a subgraph of radius $k$, this implies that $d(v_m,v) = \lambda$ for all $m \geq k$, and so $T' \in O(T)$. Therefore, we have precisely two orbits in $M$, $O(T)$ and $O(\cG_n)$.

The orbit $O(\cG_n)$ is periodic, and the orbit $O(T)$ is countable. Therefore, $\dim_H M = 0$.

We show that $M$ has positive box dimension. Since the box dimension of a set is equal to the box dimension of its closure, it is enough to compute the box dimension of $O(T)$. Let $N = card(V(P))$, the number of vertices in $P$. 

Denote by $B_i$ the number of clopen balls of diameter $e^{-i}$ required to cover $M$. This number is equal to the number of distinct patterns of radius $i$, which occur in $T$. We now obtain the lower and the upper bounds on $B_i$.

Suppose a vertex $g$ is in $T$ but not in $P$, and $i > 2$. If $d(g,v) < i$, then the pattern $B_T(g,i)$ contains the vertex $v$, which is attached to two edges. Then for every other $g' \in T - P$ the relative location of $v$ in the pattern $B_T(g',i)$ will be different from that in $B_T(g,i)$, so all such patterns are pairwise distinct. The number of such patterns is equal to at least the number of vertices between the $\partial_{\ell} P$ and $\partial_i B_T(v,i)$, which is
   $$1 + (2n-1) + (2n-1)^2 + \cdots + (2n-1)^{i-1} =  \frac{(2n-1)^{i} - 1}{2n -2}.$$
 If $d(g,v) > i$, then  there is an isomorphism $B_{T}(g,i) \to B_{\cG_n}(e,i)$. Since the number of distinct patterns, centered at vertices, contained in $P$, is at most $N$, we obtain the following estimates on the number $B_i$ of clopen balls of diameter $e^{-i}$, required to cover $M$:
  \begin{align*} (2n-1)^{i-2} \leq \frac{(2n-1)^{i} - 1}{2n-2}+1 \leq B_i \leq N+  \frac{(2n-1)^{i} - 1}{2n-2} +1 \leq (2n-1)^{i+\alpha} \end{align*}  
for some $\alpha >1$,  which implies that
  \begin{align*} \underline{\dim}_B M = \overline{\dim}_B M = \ln (2n-1).\end{align*}
\endproof

 The following result shows statement $3)$ of Proposition \ref{prop-subsets}.

\begin{prop}\label{positive}
The set $\cJ$ of closed invariant transitive sets with non-zero Hausdorff dimension is dense in $X_n$.
\end{prop}

\proof Let $\cA$ be a finite set, and let $\Sigma = \{\sigma: \mZ \to \cA\}$ be a shift space. Let $\Phi: \Sigma \to X_n$ be the embedding defined by \eqref{eq-embedmap} and \eqref{eq-themap}, and denote $M_1 = \Phi(\Sigma)$. Let $\sigma \in \Sigma$ be an element with dense orbit in $\Sigma$, and let $(T_1,e) = \Phi(\sigma) \in M_1$. Then the closure $\overline{O(T_1)} = M_1$. Recall that by Lemma \ref{shift} $\dim_H(M_1) > 0$.

Now let $U$ be a clopen neighborhood in $X_n$, and let $T_2 \in U$. Let $M_2 = \overline{O(T_2)}$. In \cite{Lukina2012}, the author introduced the \emph{fusion} of graphs. Given two graphs $T_1$ and $T_2$, the fusion produces a graph $\cT \in X_n$ such that $T_1,T_2 \subset \cM = \overline{O(\cT)}$, i.e. $T_1$ and $T_2$ are both in the closure of the orbit of $\cT$. 

This is achieved as follows: let $a,b$ be any two of $n$ generators in $F_n^+$, and let $K$ be an infinite subgraph of $\cG_n$, containing the vertices $\{a^k,b^k \mid k \in \mZ\}$. That is, $K$ is a union of two lines, one of which is made up of edges marked by $a$, and another one is made up of edges marked by $b$. The graph $K$ contains the identity vertex in the intersection. 

Next, for $\ell >1$ let $P_\ell=B_{T_1}(e,2^\ell)$ and $Q_\ell = B_{T_2}(e,2^\ell)$, so that $\{P_\ell\}_{\ell \geq 1}$ is an increasing nested sequence of finite graphs which exhausts $T_1$, and $\{Q_\ell\}_{\ell \geq 1}$ is an increasing nested sequence of finite graphs which exhausts $T_2$. To fuse $T_1$ and $T_2$, we translate $\{P_\ell\}_{\ell \geq 1}$ and $\{Q_\ell\}_{\ell \geq 1}$ inside the Cayley graph $\cG_n$ in such a way that the translated graphs are pairwise disjoint and their union with $K$ is path connected. This can always be done, and a detailed construction is described in \cite{Lukina2012}. We denote the union of $K$ with the translated graphs $\{P_\ell\}_{\ell \geq 1}$ and $\{Q_\ell\}_{\ell \geq 1}$ by $\cT$ and say that we obtained $\cT$ by attaching $P_\ell$ and $Q_\ell$, $\ell \geq 1$, to $K$. Depending on the labels of edges in the boundaries of $P_\ell$, $P_\ell$ are attached either to the half line in $K$ which contains the vertices $a^n$ with $n \geq 1$, or to the half line in $K$ which contains the vertices $b^n$ with $n \geq 1$, for $\ell \geq 1$. Similarly, depending on the labels of edges in the boundaries of $Q_\ell$,  $Q_\ell$ are attached either to the half line in $K$ which contains the negative powers $a^{-n}$ as vertices, or to the half line in $K$ which contains the negative powers $b^{-n}$ as vertices, for $n \geq 1$ and $\ell \geq 1$.
The resulting graph $\cT$ has a property that it contains all patterns around the distinguished vertex which occur in $T_1$, and also all patterns around the distinguished vertex which occur in $T_2$ infinite number of times. It follows that $T_1,T_2 \in \overline{O(\cT)}$, and so are their orbits under the pseudogroup action.

\begin{figure}[!htbp] 
\centering
    \includegraphics[width=0.7\textwidth]{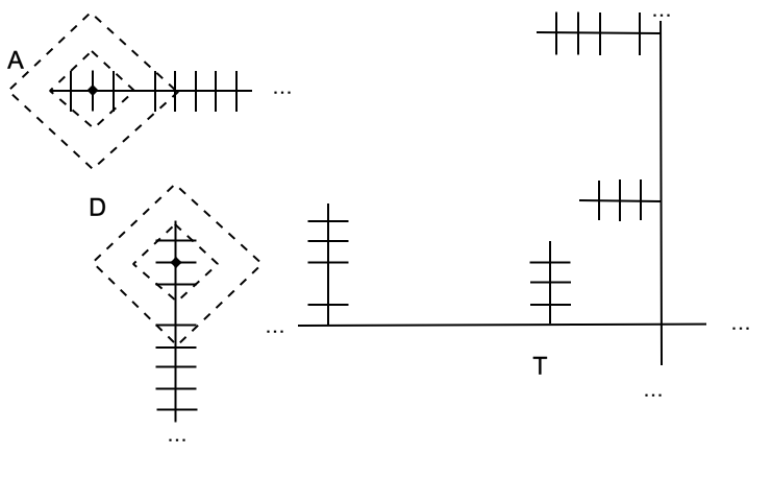}
\caption{Applying fusion to the graphs $A$ and $D$ to obtain the tree $T$. For simplicity we chose subgraphs of $\cG_2$ which can be embedded into the plane in such a way that the embedded edges have equal length.}    
\label{fig:fusion}
\end{figure}

The construction of fusion is illustrated in Figure \ref{fig:fusion} for the case $n=2$. For simplicity we chose subtrees in $\cG_2$ which can be embedded into the plane so that the edges of the embedded graphs have equal lengths (an embedding with equal edge length property need not exist for an arbitrary subtree of $\cG_2$). Dashed lines in Figure \ref{fig:fusion} correspond to the boundaries of $P_\ell$ and $Q_\ell$, $\ell \geq 1$. Their translated copies are attached to the graph $K$ to obtain the tree $T$.

Since $O(T_1)$ is dense in $M_1$, then $M_1\subset \cM$, and by \cite[Theorem 6.1]{Pesin1997} $\dim_H \cM\geq \dim_H M_1 >0$. Since $T_2 \in U$, then $\cM \cap U \ne \emptyset$. Thus the set of closed invariant transitive sets with non-zero Hausdorff dimension is dense in $X_n$.
\endproof

\section{Application: graph matchbox manifolds} \label{sec-application}

In this section, we prove Proposition \ref{prop-main3}  and give a brief overview of the literature on the embeddings of laminations as subsets of foliations. 

Since we are interested in embeddings into differentiable foliations, we start by outlining a few properties which such embeddings must possess. For details on foliated manifolds, the reader can refer to the texts by 
Moore and  Schochet \cite{MS2006} or by Candel and Conlon \cite{CandelConlon2000}, and for the properties of embeddings into differentiable foliations, to the paper by Hurder \cite{Hurder2013}.

\medskip
Let $M$ be a compact $C^\infty$ manifold of dimension $m$, and let $\cF$ be a $C^{\infty,r}$-foliation of $M$, where $r \geq 1$ is an integer. 
Here $C^{\infty,r}$ means that the leaves of a foliation are smooth manifolds, but the manifold $M$ admits a foliated atlas of class $C^r$, $r \geq 1$. The requirement that $r \geq 1$, that is, the transition maps of any foliated atlas on $M$ are differentiable, imposes important restrictions on the dynamics of the foliation.

A smooth manifold $M$ has a Riemannian metric, which induces a metric $d_M$ on $M$. By a standard result in foliation theory, one can choose a finite foliated atlas
 $$\cU = \{\varphi_i: U_i \to (-1,1)^{m-q} \times V_i\}_{1 \leq i \leq \alpha}, ~ V_i \subset \mR^q \textrm{ is open },$$ 
 of $M$ which is \emph{regular}: that is, the topological closure $\overline{U}_i$ of each chart  is contained in an open metric ball around the center $\varphi_i^{-1}(0,0)$ of the chart, the open set $U_i$ is equal to the interior of $\overline{U}_i$, connected components of leaves in $U_i$ are geodesically convex sets which intersect in simply connected sets. Also, without loss of generality we can assume  that the local transversals ${\ds T_i = \varphi^{-1}_i(\{0\} \times (-1,1)^q)}$ are disjoint. The Riemannian metric on $M$ restricts to a Riemannian metric on each $T_i$, $1 \leq i \leq \alpha$, and induces a distance function on $T_i$. Pushing forward this distance function via the corresponding homeomorphism $\varphi_i$, $1 \leq i \leq \alpha$, to the transversal $V_i$, we obtain a distance function $d_{V_i}$ on $V_i$. We now can define a metric $d_V$ on $V$  by saying that $d_V(u,v) =  d_{V_i}(u,v)$ if $u,v \in V_i$, and $d_V(u,v) = \max\{{\rm diam}(V_k) \mid 1 \leq k \leq \alpha\} +1$ if $u \in V_i$ and $v \in V_j$.

If two charts $U_i$ and $U_j$ in the atlas have a non-empty intersection, $U_i \cap U_j \ne \emptyset$, then there is a $C^r$-diffeomorphism $g_{ij}$ from an open subset of $V_i$ to an open subset of $V_j$. We call the collection $\{g_{ij}\}_{1 \leq i,j \leq \alpha}$, $r \geq 1$, of $C^r$-diffeomorphisms corresponding to all non-empty intersections of charts the \emph{generating set} of the holonomy pseudogroup $\G$ of $\cF$, and the elements in the collection the \emph{generators}. The holonomy pseudogroup $\G$ contains all possible compositions of the generators, the restrictions of compositions to open subsets of their domains, and glueings of two compositions into one map, if these compositions coincide on the intersection of their domains.

Since the holonomy group is generated by $C^r$-diffeomorphisms, for $r \geq 1$, the holonomy pseudogroup $\G$ is $C$-\emph{Lipschitz} with respect to $d_V$ (see Ghys, Langevin and Walczak \cite{GLW1988} or Hurder \cite[Proposition 4.2]{Hurder2013} for a proof). This means that there exists a constant $C>0$ such that for every generator $\{g_{ij}\}_{1 \leq i,j \leq \alpha}$ and every $x,y \in {\rm dom}( g_{ij})$ one has 
\begin{align*} \frac{1}{C} d_{V}(x,y) \leq d_V(g_{ij}(x),g_{ij}(y)) \leq C d_V(x,y) . \end{align*}
In other words, each generator $g_{ij}$ is a \emph{bi-Lipschitz} map, with a common constant $C>0$.

Now suppose that $\cS$ is a closed saturated subset of $M$, that is, $\cS$ is closed in the topology of $M$, and if a leaf $L$ of the foliation of $M$ intersects $\cS$ non-trivially, then $L \subseteq \cS$. Then every $C^{\infty,r}$ atlas $\cU$ of $M$, for $r\geq 1$, induces a $C^{\infty,0}$ atlas $\cU' = \{\varphi_i: U_i \cap \cS \to [-1,1]^{m-q} \times Z_i\}$ of $\cS$, where $Z_i$ is in general a topological space and need not carry any transverse differentiable structure. 

The metric $d_M$ induces a metric on $\cS$ by restriction. Since each $Z_i$ is a subset of $V_i$, the metric $d_V$ induces a metric $d_Z$ on $Z=Z_1 \cup \cdots \cup Z_\alpha$ by restriction, and the holonomy pseudogroup $\G_{\cS}$ is $C$-Lipschitz relative to $d_Z$. Therefore, we have the following criteria.

If a lamination $\fM$ can be realized as a subset of a $C^{\infty,r}$ foliation of a manifold $M$, for $r \geq 1$, then it must satisfy the following two conditions:
\begin{enumerate}
\item $\fM$ must admit a foliated atlas $\cQ$ and a metric $d_{\fM}$, such that the holonomy pseudogroup associated to $\cQ$ is $C$-Lipschitz with respect to the metric induced on the transversal by $d_{\fM}$, for some constant $C>0$.
\item There must exist a homeomorphism onto its image $\phi: \fM \to M$ which is a bi-Lipschitz map, that is, there exists a constant $K > 0$ such that for any $x,y \in \fM$ we have
\begin{align}\label{eq-lip} \frac{1}{K} d_{\fM}(x,y) \leq d_M(\varphi(x),\varphi(y)) \leq K d_\fM(x,y) . \end{align}
If $T$ and $Z$ are transversals for $M$ and $\fM$ respectively, then \eqref{eq-lip} implies the estimates
\begin{align}\label{eq-liptr} \frac{1}{K'} d_{Z}(x,y) \leq d_T(\varphi(x),\varphi(y)) \leq K' d_Z(x,y), \end{align}
for some $K' >0$.
\end{enumerate}

\medskip
Now consider the space of pointed trees $X_n$, associated to the Cayley graph of a free group $F_n$, $n \geq 2$, with the ball metric $d_n$.
Recall \cite{Ghys1999,ALM2009a,Blanc2001,LozanoLukina2012,Lukina2012} that the pseudogroup dynamical system $(X_n,d_n,\G_n)$, $n \geq 2$, can be suspended to produce a lamination $\fM_n$ with leaves of dimension $2$ as in the following theorem. 

\begin{thm}\cite{Ghys1999}\label{thm-genconstruction}
Let $F_n$ be a free group on $n \geq 2$ generators, and $(X_n,d_n)$ be the corresponding space of pointed trees with the action of a pseudogroup $\G_n$. Then there exists a compact metric space $\fM_n$, and a finite smooth foliated atlas $\cV = \{\phi_i:V_i \to U_i \times \fX_i\}_{1 \leq i \leq \nu}$, where $U_i \subset \mR^2$ is open, and $\fX_i$ is a Cantor set, with associated holonomy pseudogroup $\fP$, such that the following holds. 
\begin{enumerate}
\item The leaves of $\fM_n$ are Riemann surfaces.
\item There is a homeomorphism onto its image
  \begin{align*} t: X_n \to \bigcup_{1 \leq i \leq \nu} \fX_i, \end{align*}
such that $t(X_n)$ is a transverse model for the foliation $\cF$, and $\fP|_{\tau(X_n)} = t_*\G$, where $t_*\G$ is the pseudogroup induced on $t(X_n)$ by $\G$.
\end{enumerate}
\end{thm}

The proof of Theorem \ref{thm-genconstruction} can be found in \cite{Ghys1999,Blanc2001}. We remark that by a small modification of the proof of Theorem \ref{thm-genconstruction} one can obtain a suspension $\fM_n$ of the pseudogroup dynamical system $(X_n,d_n,\G_n)$ with leaves of dimension $3$ and higher. We do not discuss the details of the construction in the proof of Theorem \ref{thm-genconstruction} here since the understanding of our results only requires the knowledge of the transverse dynamics of a lamination.

A \emph{graph matchbox manifold} is a compact connected subspace $\cM$ of the lamination $\fM_n$ given by Theorem \ref{thm-genconstruction}, which is a closure of a leaf $L$ in $\fM_n$, i.e. $\cM = \overline{L}$. 
 It was shown in \cite{Lukina2012} that $\fM_n$ contains a residual subset of dense leaves, and so $\fM_n$ is a graph matchbox manifold in its own right.  Then 
Proposition \ref{prop-subsets} leads to the following non-embedding result.

\begin{thm}
Let $F_n$, $n \geq 2$, be a free group on $n$ generators, and $(X_n,d_n)$ be the associated space of pointed graphs with ball metric. Let $\fM_n$ be a suspension of the natural pseudogroup action on $(X_n,d_n)$ as in Theorem \ref{thm-genconstruction}. Then $\fM_n$ cannot be embedded as a subset of a $C^\infty$-manifold with a $C^{\infty,r}$-foliation by a bi-Lipschitz homeomorphism onto its image, for $r \geq 1$.
\end{thm}

\proof Let $M$ be a foliated manifold with foliation $\cF$, and let $Z$ be any transverse $k$-dimensional section. By Statement (2) in \cite[Theorem 6.1]{Pesin1997} the Hausdorff dimension of any subset of $Z$ is at most $k$. 
Suppose there exists a foliated bi-Lipschitz embedding $\phi: \fM_n \to M$, then it must map $X_n$ onto a transverse section $Z$, or a finite union of such sections. By Proposition \ref{prop-subsets} the Hausdorff dimension of $X_n$ is infinite, so the Hausdorff dimension of the image $\phi(X_n) \subset Z$ must be infinite, a contradiction. Therefore, no such bi-Lipschitz homeomorphism $\phi$ can exist.
\endproof

To author's best knowledge, this example is the first example of a lamination with a given metric and a $C$-Lipschitz pseudogroup (see Remark \ref{rem-lipschitz} for a proof that $\G_n$ is $C$-Lipschitz), which does not admit a bi-Lipschitz embedding into a $C^{\infty,r}$-foliation of a finite-dimensional manifold, for $r \geq 1$.

Indeed, very few examples in the modern literature on foliations deal with embeddings of laminations into differentiable foliations. Most of the available results are about continuous embeddings of various types of laminations into Euclidean spaces and into $C^{\infty,0}$-foliations. The properties of these embeddings are very different from the embeddings into $C^{\infty,r}$ foliations, for $r \geq 1$. 

For embeddings into Euclidean spaces, the embedding map is a homeomorphism onto its image, but it need not be bi-Lipschitz. For embeddings into $C^{\infty,0}$-foliations, the embedding map is a homeomorphism onto its image, not necessarily bi-Lipschitz, and, additionally, the holonomy homeomorphisms on the embedded Cantor set transversal must extend to homeomorphisms of a transversal $T$ of the foliated manifold $M$. The extensions of homeomorphisms are not required to be differentiable.

Embeddings of solenoids into Euclidean spaces were studied in continuum theory. Recall \cite{McCord1965} that a
\emph{regular} solenoid is a compact metrizable space, homeomorphic to the inverse limit of a sequence of finite-to-one coverings of a closed manifold $M_0$, such that the fundamental group of each covering space $M_i$ in the sequence injects onto a proper \emph{normal} subgroup of $M_0$. Regular solenoids 
are the simplest examples in the class of weak solenoids, which are the spaces homeomorphic to inverse limits of sequences of any finite-to-one coverings of a closed manifold. Regular solenoids were introduced by McCord \cite{McCord1965}, with details and examples of regular and weak solenoids given in Schori \cite{Schori1966}, Rogers and Tollefson \cite{RogersTollefsonf1971}, Fokkink and Oversteegen \cite{FO2002}. Every solenoid is a fibre bundle with Cantor set fibre over any closed manifold in the corresponding sequence, and so every solenoid is an example of a lamination.

The study of embeddability of solenoids into manifolds dates back to Bing \cite{Bing1960} who showed that the inverse limit of non-trivial self-coverings of a circle cannot be topologically embedded into a plane. 
Prajs \cite{Prajs1990} showed that, if a homogeneous solenoid with leaves of dimension $n$ is embedded into an $(n+1)$-dimensional manifold, then it is an $n$-dimensional manifold itself. Rephrasing, a non-trivial homogeneous solenoid cannot be embedded into a manifold with codimension $1$.
Clark and Fokkink \cite{ClarkFokkink2004} proved a similar result using a cohomological argument. Jiang, Wang and Zheng \cite{JWZ2008} extended the non-embedding results of \cite{ClarkFokkink2004,Prajs1990} to all weak solenoids.

Keesling and Wilson \cite{KeeslingWilson1985} proved, that inverse limits of coverings of $n$-tori embed into $\mathbb{R}^{n+2}$.
Clark and Fokkink \cite{ClarkFokkink2004} gave necessary and sufficient conditions for a regular solenoid over a closed manifold $M_0$ to be topologically embeddable into $M_0 \times \mathbb{R}^2$. 
In the same paper, the authors also considered topological embeddings of regular solenoids into foliated bundles. The holonomy action on the fibre of a solenoid is a group action. Under certain conditions one can obtain an embedding of a regular solenoid into $M_0 \times \mathbb{R}^2$, so that the product structure on the target space and the bundle structure match. By the result of Moise \cite{Moise1977} any homeomorphism of a Cantor subset of $\mR^2$ can be extended to a homeomorphism of $\mR^2$. Using that, to each generator of the fundamental group of $M_0$ one can associate a homeomorphism of $\mR^2$ to obtain a group $G \subset {\rm Homeo}(\mR^2)$. The homeomorphisms, obtained by such extension, need not satisfy the same relations as the generators of the fundamental group of $M_0$, except when restricted to the embedded fibre. Thus the group $G$ acting on the fibre $\mR^2$ of the product $M_0 \times \mR^2$ is a free group. The embedding of a solenoid, obtained by this method, is $C^{\infty,0}$.

Lozano \cite{Lozano2009} used a similar construction, based on the result of Moise, to construct codimension $2$ embeddings of graph matchbox manifolds into a foliated manifold. The embeddings in \cite{Lozano2009} are $C^{\infty,0}$. Hector \cite{Hector2014} states that certain minimal graph matchbox manifolds without holonomy cannot embed into $C^{\infty,0}$ foliations of codimension $1$. Our result in Proposition \ref{prop-main3}  shows, that $\fM_n$ cannot be embedded into a $C^{\infty,r}$-foliation of a smooth manifold for $r \geq 1$ and for \emph{any} codimension, thus showing that the embeddings of graph matchbox manifolds into differentiable foliations are a lot more subtle than their embeddings into continuous foliations.

Williams \cite{Williams1967,Williams1974} showed that shift maps of generalized solenoids model topologically expanding hyperbolic attractors, and Brown \cite{Brown2010} showed that toral solenoids provide topological models for non-expanding hyperbolic attractors. Clark and Hurder \cite{ClarkHurder2010a} constructed embeddings of a certain class of toral solenoids as subsets of $C^{\infty,r}$-foliated manifolds, for $r \geq 1$. 

Further examples of laminations arise as the tiling spaces of aperiodic tilings of $\mathbb{R}^n$ with finite local complexity, of which there is a large literature. See, for instance, Sadun \cite{Sadun2008} and Kwapisz \cite{Kwapisz2011}.

There has been a few publications about the existence of bi-Lipschitz embeddings of tiling spaces into Euclidean spaces, which use Hausdorff and Assouad dimensions as obstructions to the existence of such embeddings. Those are papers by Julien and Savinien \cite{JulienSavinien2011} and by Bellissard and Julien \cite{BJ2012}. See also references in \cite{BJ2012} for the literature on bi-Lipschitz embeddings of ultrametric spaces into Euclidean spaces, and \cite{Fraser2020} for the Assouad dimension.

Julien and Savinien \cite{JulienSavinien2011} showed that a certain class of tiling spaces of tilings of $\mathbb{R}^n$ can be embedded into $\mathbb{R}^{n+1}$ by a bi-Lipschitz homeomorphism. They also computed the Hausdorff dimension of the canonical transversal for certain classes of tiling spaces.
Julien and Bellissard \cite{BJ2012} concentrate on embeddings of ultrametric Cantor sets, and show that linearly repetitive subshifts, Sturmian subshifts and some other examples can be embedded into a Euclidean space by a bi-Lipschitz homeomorphism. 

Hurder \cite{Hurder2013} studied Lipschitz properties of holonomy pseudogroups from the point of view of embeddings into $C^{\infty,r}$-foliations, for $r \geq 1$. We recommend his article for a nice overview of the topic.

\section{Acknowledgements} 

The author thanks Steve Hurder for useful discussions. The author is thankful to the anonymous referee for the careful reading of the manuscript and comments which improved the exposition.


\end{document}